\documentclass[12pt]{article}
\usepackage{amsmath,amsthm,verbatim,amssymb,amsfonts,amscd,diagrams, graphics, mathrsfs}
\theoremstyle{plain}
\newtheorem{theorem}{Theorem}[section]
\newtheorem{corollary}[theorem]{Corollary}
\newtheorem{lemma}[theorem]{Lemma}
\newtheorem{proposition}[theorem]{Proposition}

\theoremstyle{definition}
\newtheorem{definition}[theorem]{Definition}

\theoremstyle{remark}
\newtheorem{remark}[theorem]{Remark}

\newarrow{ul}---->
\newarrow{Backwards}<----

\newcommand{\td}[1]{\tilde{#1}}
\newcommand{\into}{\hookrightarrow}
\newcommand{\Z}{\mathbb{Z}}
\newcommand{\Q}{\mathbb{Q}}
\newcommand{\R}{\mathbb{R}}

\newcommand{\F}{\mathbb{F}}
\newcommand{\D}{\mathbb{D}}

\newcommand{\bd}{\partial}

\renewcommand{\H}{\mathbb H}

\newcommand{\mc}[1]{\mathcal{#1}}
\newcommand{\ms}[1]{\mathscr{#1}}

\newcommand{\dlim}{\varinjlim}

\newcommand{\hl}{\operatorname{holink}}

\newcommand{\Hom}{\text{Hom}}
\newcommand{\Ext}{\text{Ext}}

\newcommand{\im}{\text{im}}

\newarrow{Onto}----{>>}
\newarrow{Equals}=====

\begin{document}

\title{Intersection homology and Poincar\'e duality on homotopically stratified spaces}
\author{Greg Friedman\\Texas Christian University}
\date{April 17, 2007}
\maketitle

\begin{abstract}
We show that intersection homology extends Poincar\'e duality to manifold homotopically stratified spaces (satisfying mild restrictions). This includes showing that, on such spaces, the sheaf of singular intersection chains is quasi-isomorphic to the Deligne sheaf. 
\end{abstract}


\tableofcontents

\section{Introduction}
The primary purpose of this paper is to extend Poincar\'e duality to manifold homotopically stratified spaces using intersection homology.

Intersection homology was introduced by Goresky and MacPherson \cite{GM1} in order to extend Poincar\'e duality to \emph{manifold stratified spaces} -- spaces that are not manifolds but that are composed of manifolds of various dimensions. Initially, this was done for \emph{piecewise-linear pseudomanifolds} \cite{GM1}, which include algebraic and analytic varieties\footnote{excluding those with codimension one strata}, but the result
was soon extended to \emph{topological pseudomanifolds} (Goresky-MacPherson \cite{GM2}) and \emph{locally conelike topological stratified spaces}, also called \emph{cs-spaces} (Habegger-Saper \cite{HS91}).  The establishment of  Poincar\'e duality for pseudomanifolds has led to the successful study and application  of further related invariants. To name just a few: Right in \cite{GM1}, Goresky and MacPherson introduced signatures and $L$-classes for pseudomanifolds with only even codimension strata;  Siegel extended signatures and bordism theory to Witt spaces \cite{Si83}; and various extensions of duality and characteristic classes have been studied by  Cappell, Shaneson, and Banagl, in various combinations  \cite{CS, Ba02, Ba06, BCS}.  For applications of interesection homology in this direction, we refer the reader to \cite{BaIH}; for applications of intersection homology in other fields, we refer the reader to \cite{KIR}.

In \cite{Q1}, Quinn introduced   \emph{manifold homotopically stratified spaces} (MHSSs), with the intent to provide ``a setting for the study of purely topological stratified phenomena, particularly group actions on manifolds.'' 
In this context of topological group actions on manfiolds, MHSSs have been  studied by Yan \cite{Ya93}, Beshears \cite{Be97}, and  Weinberger and Yan \cite{Ya01, Ya05}.\footnote{
The application of intersection homology to the  study of group actions  both on smooth manifolds  and on stratified spaces is an active field of research; see, e.g.,  \cite{HS93, Sa96, Bry, Cu92, PS04, Pa05}.} But MHSSs also arise in categories with more stucture; for example, Cappell and Shaneson showed that they occur as  mapping cylinders of maps between smoothly stratified spaces \cite{CS95}. A surgery theory for MHSSs has been developed by Weinberger \cite{Wein}, and their  geometric neighborhood properties have been studied by Hughes, culminating in \cite{Hu02}.
In \cite{Q2}, Quinn noted that MHSSs ``are defined by local homotopy properties, which seem more appropriate for the study of a homology theory'' than the local homeomorphism properties of pseudomanifolds, and he showed that intersection homology is a topological invariant on these spaces, independent of the stratification. A further survey of MHSSs can be found in Hughes and Weinberger \cite{HW}.

We prove the following Poincar\'e duality theorem, which generalizes the Goresky-Siegel extension of Goresky and MacPherson's intersection homology duality. The stated neighborhood condition on the MHSS $X$ is described more fully below in Section \ref{S: LATN} but includes MHSSs with compact singular set $\Sigma$ such that all non-minimal strata of $X$ have dimension $\geq 5$. The condition of being \emph{homotopy locally ($\bar p$,$R$)-torsion free} is a weakening of Goresky and Siegel's \emph{locally $\bar p$-torsion free} and is defined in Section \ref{S: pd}; roughly, this condition requires the torsion to vanish from certain local intersection homology groups that, for a pseudomanifold, would correspond to certain intersection homology groups of the links.

\begin{theorem}[Theorem \ref{T: PD}]
Let $X$ be a homotopy locally ($\bar p$,$R$)-torsion free n-dimensional MHSS with no codimension one stratum  and with sufficiently many local approximate tubular neighborhoods (in particular, if all non-minimal strata of $X$ have dimension $\geq 5$). Let $\mc O$ be the orientation sheaf of the $n$-manifold $X-X^{n-2}$, and let $\mc E$ be a local coefficient system on $X-X^{n-2}$ of finitely-generated free modules over the  principal ideal domain $R$. Let $\bar p$ and $\bar q$ be dual perversities ($\bar p(k)+\bar q(k)=k-2$).
Let $TH_*$ and $FH_*$ denote, respectively, the $R$-torsion subgroup and $R$-torsion free quotient group of $IH_*$, and let $Q(R)$ denote the field of fractions of $R$.

Suppose that $\Hom(T^{\bar p}H^c_{i-1}(X;\mc E),Q(R)/R)$ is a torsion $R$-module (in particular, if $T^{\bar p}H^c_{i-1}(X;\mc E)$ is finitely generated).
Then $$\Hom(F^{\bar p}H^c_i(X;\mc E),R)\cong F^{\bar q}H^{\infty}_{n-i}(X;\Hom(\mc E, R_{X-X^{n-2}})\otimes \mc O)$$ and $$\Hom(T^{\bar p}H^c_{i-1}(X;\mc E),Q(R)/R)\cong T^{\bar q}H^{\infty}_{n-i}(X;\Hom(\mc E, R_{X-X^{n-2}})\otimes \mc O).$$  
\end{theorem}

We record separately the case for field coefficients, for which all of the  torsion conditions are satisfied automatically.

\begin{corollary}[Corollary \ref{C: PD}]
Let $X$ be an n-dimensional MHSS with no codimension one stratum  and with sufficiently many local approximate tubular neighborhoods (in particular, if all non-minimal strata of $X$ have dimension $\geq 5$). Let $\mc O$ be the orientation sheaf of the $n$-manifold $X-X^{n-2}$, and let $\mc E$ be a local coefficient system on $X-X^{n-2}$ of finitely-generated $\F$-modules for a field $\F$. Let $\bar p$ and $\bar q$ be dual perversities ($\bar p(k)+\bar q(k)=k-2$). Then $$\text{\emph{Hom}}(I^{\bar p}H^c_{n-i}(X;\mc E);\F)\cong I^{\bar q}H^{\infty}_i(X;\Hom(\mc E, \F_{X-X^{n-2}})\otimes \mc O).$$
\end{corollary}

In particular, if $X$ is a compact orientable MHSS satisfying the hypotheses of the corollary, we obtain the more familiar pairing $$\Hom(I^{\bar p}H_i(X;\Q),\Q)\cong I^{\bar q}H_{n-i}(X;\Q).$$ If, in addition, $X$ is homotopy locally ($\bar p$,$\Z$)-torsion free, we have
$$\Hom(F^{\bar p}H_i(X),\Z)\cong F^{\bar q}H_{n-i}(X) \qquad \text{and} \qquad \Hom(T^{\bar p}H_i(X),\Q/\Z)\cong T^{\bar q}H_{n-i}(X),$$
which generalize the usual intersection and linking pairings for manifolds.

In the final section of the paper, Section \ref{S: rings}, we explore conditions that would ensure an appropriate duality over more general ground rings. 

\begin{comment}
\begin{theorem}[Theorem \ref{T: PD}]
Let $X$ be an n-dimensional MHSS with no codimension one stratum  and with sufficiently many local approximate tubular neighborhoods. Let $\mc O$ be the orientation sheaf of the $n$-manifold $X-X^{n-2}$, and let $\mc E$ be a local coefficient system of $\F$-modules for a field $\F$. Let $\bar p$ and $\bar q$ be dual perversities ($\bar p(k)+\bar q(k)=k-2$). Then $$I^{\bar q}H^{\infty}_i(X;\Hom(\mc E, \F_{X-X^{n-2}})\otimes \mc O)\cong \text{\emph{Hom}}(I^{\bar p}H^c_{n-i}(X;\mc E);\F).$$
\end{theorem}

In particular, if $X$ is a compact orientable MHSS satisfying the hypotheses of the theorem, we obtain nonsingular pairings 
$$I^{\bar p}H_i(X;\Q)\otimes I^{\bar q}H_{n-i}(X;\Q)\to \Q.$$

\end{comment}

\paragraph{Outline.} Section \ref{S: philosophy} contains a brief overview of our need to utilize  singular chain intersection homology on MHSSs and the relationship between this approach and the sheaf-theoretic point of view. 
In Section \ref{S: bg}, we review the requisite background technical material. In Section \ref{S: LATN}, we establish the notion of MHSSs with sufficiently many approximate tubular neighborhoods; these are the spaces on which our results will apply. Section \ref{S: deligne} contains the proof that the sheaf complex of singular intersection chains on an MHSS is quasi-isomorphic to the Deligne sheaf of \cite{GM2}. In Section \ref{S: cc}, we demonstrate that these sheaves are constructible on MHSSs, and in Section \ref{S: pd}, we establish Poincar\'e duality. Section \ref{S: apps} contains a definition of Witt spaces in the class of MHSSs, and Section \ref{S: rings} explores how our duality results extend for ground rings of higher cohomological dimension. Finally, we provide a technical computation in the Appendix.

\paragraph{Acknowledgments.} I thank Bruce Hughes and Markus Banagl each for several helpful discussions.

\section{Sheaves vs. Singular Chains}\label{S: philosophy}

 Intersection homology on piecewise linear (PL) pseudomanifolds was defined  initially in terms of simplicial chains, and the original proof of Poincar\'e duality over a field on compact orientable PL pseudomanifolds in \cite{GM1} was performed via a combinatorial construction using the triangulations of the spaces. However, by \cite{GM2} sheaf theory had taken over. It was shown that intersection homology on a PL pseudomanifold can be obtained as the hypercohomology of a certain sheaf complex, the Deligne sheaf complex, that can be defined without any reference to simplicial chains. This sheaf theoretic description of intersection homology was extremely successful - by eliminating a need for simplicial chains, intersection homology could be extended to topological pseudomanifolds (which need not be triangulable), and the existence of an axiomatic characterization of the Deligne sheaf led to purely sheaf-theoretic proofs of topological invariance (independence of stratification) and of  Poincar\'e duality, via Verdier duality of sheaf complexes. It is in this sheaf-theoretic realm that many of the most important applications of intersection homology have been attained, including intersection homology versions of Hodge theory, the Lefschetz hyperplane theorem, and the hard Lefschetz theorem for singular varieties, as well as applications to the Weil conjecture for singular varieties, the Riemann-Hilbert correspondence, and the proof of the Kazhdan-Lusztig conjecture  (see \cite{KIR} for an exposition of these applications). 

However, the homotopy theoretic nature of MHSSs makes it difficult to work with sheaves on these spaces - sheaf cohomology does not always behave well with respect to homotopies and homotopy equivalences, and some of the spaces that occur in the analysis of MHSSs, such as certain path spaces, are not locally-compact, a property that is often required in order to employ some of the most useful theorems of sheaf theory. This discourages one from taking a purely sheaf theoretic approach to intersection homology on these spaces. Fortunately, a singular chain version of intersection homology exists, due initially to King \cite{Ki}, and this is the version of intersection homology that  Quinn demonstrated was a topological invariant of MHSSs in \cite{Q2}. In \cite{GBF10}, we showed that the singular intersection chiains also generate a sheaf complex and that it is quasi-isomorphic to the Deligne sheaf  complex on topological pseudomanifolds. In this paper, we will show that that these two versions of intersection homology also agree on MHSSs, and then we will \emph{use} this correspondence to obtain Poincar\'e duality by applying chain-theoretic arguments to demonstrate properties of the sheaf complexes.

More specifically, the Goresky-MacPherson proof of Poincar\'e duality on topological pseudomanifolds \cite{GM2} proceeds by establishing that the Deligne sheaf complex is characterized by its axiomatic properties and then by showing that the Verdier dual of a perversity $\bar p$ Deligne sheaf complex satisfies the axioms to be the Deligne sheaf complex with the complementary perversity. Duality of intersection homology then follows by general sheaf theory.  Showing that the Verdier dual of a Deligne sheaf satisfies the axioms to be another  Deligne sheaf is done purely sheaf-theoretically and relies at certain points upon the geometric form of local neighborhoods in pseudomanifolds - each point has \emph{distinguished neighborhoods} homeomorphic to $\R^{n-k}\times cL$, where $cL$ is the open cone on a lower-dimensional pseudomanifold.\footnote{In particular, the proof of duality in \cite{Bo} (which clarifies the proof in \cite{GM2}) shows that the Verdier dual of the Deligne sheaf satisifies the axioms AX2, which include a condition that this dual be stratified cohomologically locally constant ($X$-clc). The satisfaction of this condition follows from the Deligne sheaf itself being $X$-clc, the proof of which, given for pseudomanifolds in \cite{Bo}, relies on the distinguished neighborhoods. In the context of MHSSs, it is not clear that it should follow from these methods that the Deligne sheaf is $X$-clc and an alternative approach is therefore necessary; see Section \ref{S: cc}, below.} See \cite{GM2} or \cite{Bo} for details. The key difficulty for MHSSs is that these distinguished neighborhoods no longer necessarily exist, and so the local sheaf arguments of \cite{GM2} and \cite{Bo} no longer apply. To overcome this difficulty, we establish a quasi-isomorphism between the Deligne sheaf and the sheaf of singular intersection chains (of the same perversity), and then we compute locally using the singular chains, which are much better adapted to being manipulated by homotopy properties. Instead of distinguished neighborhoods, we use local versions of the Approximate Tubular Neighborhoods of Hughes \cite{Hu02}, which exist on all MHSSs satisfying mild dimension requirements (see Section \ref{S: LATN}, below).

\section{Background and Basic Terminology}\label{S: bg}

\subsection{Intersection homology}

In this section, we provide a quick review of the definition of intersection homology. For more details, the reader is urged to consult King \cite{Ki} and the author \cite{GBF10} for singular intersection homology and the original papers of Goresky and MacPherson \cite{GM1,GM2} and the book of Borel \cite{Bo} for the simplicial and sheaf definitions. Singular chain intersection homology theory was introduced in \cite{Ki} with finite chains (compact supports) and generalized in \cite{GBF10} to include locally-finite but infinite chains (closed supports). 

We recall that singular  intersection homology  is defined on any filtered space

$$X=X^n\supset X^{n-1}\supset \cdots \supset X^0\supset X^{-1}=\emptyset.$$
In general, the superscript ``dimensions'' are simply labels and do not necessarily reflect any geometric notions of dimension. We refer to $n$ as the \emph{filtered dimension} of $X$, or simply as the ``dimension'' when no confusion should arise. 
The set $X^i$ is called the $i$th \emph{skeleton} of $X$, and $X_i=X^i-X^{i-1}$ is the $i$th \emph{stratum}.

\begin{remark}\label{R: filter}
Our definition of a filtered space is more specific than that found in, e.g., \cite{Hug, Hu02} in that we require $X$ to have a finite number of strata and that the strata be totally ordered. If the skeleta $X^i$ are closed in $X$, then these spaces will also be ``stratified spaces satisfying the frontier condition'' - see \cite{Hug}.  
\end{remark}

A \emph{perversity} $\bar p$ is a function $\bar p: \Z^{\geq 1}\to \Z$ such that $\bar p(k)\leq \bar p(k+1)\leq \bar p(k)+1$. A \emph{traditional perversity} also satisfies $\bar p (1)=\bar p(2)=0$. One generally must restrict to traditional perversities in order to obtain the most important topological invariance and Poincar\'e duality results for intersection homology (see \cite{GM2, Bo, Ki, Q2}), athough many interesting results are now also known for superperversities, which satisfy $\bar p(2)>0$ (see \cite{CS, HS91, GBF10, GBF11, Sa05}).

Given $\bar p$ and $X$, one defines $I^{\bar p}C^c_*(X)$ as a subcomplex of $C^c_*(X)$, the complex of compactly supported singular chains\footnote{This is the usual chain complex consisting of finite linear combination of singular simplices, but we emphasize the compact supports in the notation to distinguish $C^c_*(X)$ from $C^{\infty}_*(X)$, which we shall also use.} on $X$, as follows: A simplex $\sigma:\Delta^i\to X$ in $C^c_i(X)$ is \emph{allowable} if $$\sigma^{-1}(X^{n-k}-X^{n-k-1})\subset \{i-k+\bar p(k) \text{ skeleton of } \Delta^i\}.$$ The chain $\xi\in C^c_i(X)$ is allowable if each simplex in $\xi$ and $\bd \xi$ is allowable. $I^{\bar p}C_*^c(X)$ is the complex of allowable chains. $I^{\bar p}C_*^{\infty}(X)$ is defined similarly as the complex of allowable chains in $C_*^{\infty}(X)$, the complex of locally-finite singular chains. Chains in $C_*^{\infty}(X)$ may be composed of an infinite number of simplices (with their coefficients), but for each such chain $\xi$, each point in $X$ must have a neighborhood that intersects only a finite number of simplices (with non-zero coefficients) in $\xi$.  $I^{\bar p}C_*^{\infty}(X)$ is refered to as the complex of intersection chains with \emph{closed supports}, or sometimes as \emph{Borel-Moore} intersection chains. See \cite{GBF10} for more details. 

The associated homology theories are denoted $I^{\bar p}H^c_*(X)$ and $I^{\bar p}H^{\infty}_*(X)$. We will sometimes omit the decorations $c$ or $\infty$ if these theories are equivalent, e.g. if $X$ is compact. We will also often omit explicit reference to $\bar p$ below, for  results that hold for any fixed perversity. 

Relative intersection homology is defined similarly, though we note that 
\begin{enumerate}
\item the filtration on the subspace will always be that inherited from the larger space by restriction, and
\item  in the closed support case, all chains are required to be locally-finite in the larger space. 
\end{enumerate}

If $(X,A)$ is such a filtered space pair, we use the notation $IC_*^{\infty}(A_X)$ to denote the allowable singular chains supported in $A$ that are locally-finite in $X$. The homology of this complex is $IH_*^{\infty}(A_X)$. Note that in the compact support case, the local-finiteness condition is satisfied automatically so we do not need this notation and may unambiguously refer to $IH_*^c(A)$. The injection  $0\to IC_*^{\infty}(A_X)\to IC_*^{\infty}(X)$ yields a quotient complex $IC_*^{\infty}(X,A)$ and a long exact sequence of intersection homology groups $\to IH_i^{\infty}(A_X)\to IH_i^{\infty}(X)\to IH_i^{\infty}(X,A)\to$. 

If $X$ and $Y$ are two filtered spaces, we call a map $f:X\to Y$
\emph{filtered} if the image of each component of a stratum of $X$
lies in a stratum of $Y$. \emph{N.B. This property is often referred to as ``stratum-preserving'', e.g. in \cite{Q1} and \cite{GBF3}. However, we must reserve  the term ``stratum-preserving'' for other common uses.}
In general, it is not required
that a filtered map take  strata of $X$ to strata of $Y$ of the same (co)dimension. However, if $f$ preserves codimension, or if  $X$ and $Y$ have the same filtered dimension and  $f(X_i)\subset Y_i$, then $f$ will induce a well-defined map on intersection homology (see \cite[Prop. 2.1]{GBF3} for a proof). In this case, we will call $f$ \emph{well-filtered}.
We call a well-filtered map $f$ a \emph{stratum-preserving homotopy equivalence} if there is a well-filtered
map $g:Y\to X$ such that $fg$ and $gf$ are homotopic to the appropriate
identity maps by well-filtered homotopies, supposing that $X\times I$ and $Y\times I$ are given the obvious product filtrations. Stratum-preserving homotopy equivalences induce intersection homology isomorphisms \cite{GBF3}. If stratum-preserving homotopy equivalences between $X$ and $Y$ exist,
we say that $X$ and $Y$ are \emph{stratum-preserving homotopy equivalent}, $X\sim_{sphe} Y$.  

In the sequel, all maps inducing intersection homology homomorphisms will clearly be well-filtered. Hence, we will usually dispense with explicit discussion of this point.  

It is shown in \cite{GBF10} that one can construct a sheaf of intersection chains $\mc{IS^*}$ on any filtered Hausdorff space $X$ such that, if $X$ is also paracompact and of finite cohomological dimension, then the hypercohomology $\H^*(\mc{IS^*})$ is isomorphic to $IH_{n-*}^{\infty}(X)$, where $n$ is the filtered dimension of $X$. If $X$ is also locally-compact, then $\H^*_c(\mc{IS^*})\cong IH_{n-*}^c(X)$. We will use some properties of these sheaves below, but we refer the reader to \cite{GBF10} for more detailed background.

\subsubsection{A note on coefficients}

Often throughout this paper we will leave the coefficient systems tacit so as not to  overburden the notation. However, except where noted otherwise, all results hold for any of the following choices of coefficients, where $R$ is any ring with unit of finite cohomological dimension:

\begin{itemize}
\item Any constant coefficient groups or $R$-modules.

\item Any local system of coefficients of groups or $R$-modules with finitely generated stalks defined on $X-X^{n-2}$ (see \cite{GM2,Bo,GBF10}).

\item If $\bar p$ is a superperversity (i.e. $\bar p(2)>0$; see \cite{CS, GBF10, GBF11}), any stratified system of coefficients $\mc G_0$ with finitely generated stalks as defined in \cite{GBF10} such that $\mc G_0|_{X-X^{n-1}}$ is a local coefficient system of groups or $R$-modules and  $\mc G_0|_{X^{n-1}}=0$. It is shown in \cite{GBF10} that this last coefficient system allows us to recover from singular chains the superperverse sheaf intersection cohomology on pseudomanifolds.

\end{itemize}

\subsection{Stratified homotopies and fibrations}

If $X$ is a filtered space, a map $f:Z\times A\to X$ is \emph{stratum-preserving along $A$} if, for each $z\in Z$, $f(z\times A)$ lies in a single stratum of $X$. If $A=I=[0,1]$, we call $f$ a \emph{stratum-preserving homotopy}. If $f:Z\times I\to X$ is only stratum-preserving when restricted to $Z\times [0,1)$, we say $f$ is \emph{nearly stratum-preserving}.

If $X$ and $Y$ are stratified spaces, a map $p:X\to Y$ is a \emph{stratified fibration} if it admits solutions to stratified lifting problems, i.e. if given a commuting diagram of maps
\begin{equation*}
\begin{CD}
Z&@>f>>& X\\
@V \times 0 VV && @VV p V\\
Z\times I &@>F>>&Y,
\end{CD}
\end{equation*}
such that  $Z$ is any space and $F$ is a stratum-preserving homotopy, there  exists a stratum-preserving homotopy $\td F:Z\times I\to X$ such that $p\td F=F$ and $\td F|_{Z\times 0}=f$. 

 See \cite{Hug, GBF3} for more on stratified fibrations.

\subsection{Manifold homotopically stratified spaces}\label{S: def MHSSs}

Even though the above definition of intersection homology applies to very general spaces, one usually needs to limit oneself to smaller classes of spaces in order to obtain nice properties. In this paper, we focus on the manifold homotopically stratified spaces introduced by Quinn and refined by Hughes. These spaces were introduced partly with the purpose in mind of being the ``right category'' for intersection homology - see \cite{Q2}.

There is disagreement in the literature as to what to call these spaces. Quinn, himself, calls them both ``manifold homotopically stratified sets'' \cite{Q1} and ``weakly stratified sets''  \cite{Q2}. Hughes \cite{Hu02} prefers the term ``manifold stratified spaces''. We use the term \emph{manifold homotopically stratified space} (MHSS), which seems to capture both that they are stratified by manifolds and that there are additional homotopy conditions on the ``gluing''. 

To define these spaces, we need some preliminary terminology. Except where noted, we take these definitions largely from \cite{Hu02}, with slight modifications to reflect the restrictions mentioned above in Remark \ref{R: filter}.

\subsubsection{Forward tameness and homotopy links}

If $X$ is a filtered space, then $Y$ is \emph{forward tame} in $X$ if there is a neighborhood $U$ of $Y$ in $X$  and a nearly-stratum preserving deformation retraction $R:U\times I\to X$ retracting $U$ to $Y$ rel $Y$. If the deformation retraction keeps $U$ in $U$, we call $U$ a \emph{nearly stratum-preserving deformation retract neighborhood (NSDRN)}. This last definition was introduced in \cite{GBF5}

The \emph{stratified homotopy link} of $Y$ in $X$ is the space (with compact-open topology) of nearly stratum-preserving paths with their tails in $Y$ and their heads in $X-Y$: $$\hl_s(X,Y)=\{\omega\in X^I\mid \omega(0)\in Y, \omega((0,1])\subset X-Y\}.$$ The \emph{holink evaluation} map takes a path $\omega\in \hl_s(X,Y)$ to $\omega(0)$. For $x\in X_i$, the \emph{local holink}, denoted $\hl_s(X,x)$, is simply the subset of paths $\omega\in \hl_s(X,X_i)$ such that $\omega(0)=x$. Holinks inherit natural stratifications from their defining spaces: $$\hl_s(X,Y)_j=\{\omega\in\hl(X,Y)\mid \omega(1)\in X_j\}.$$

If $X$ is metric and $\delta: Y\to(0,\infty)$ is a continuous function, then $\hl_s^{\delta}(X,Y)$ is the subset of paths $\omega\in\hl_s(X,Y)$ such that $\omega(I)$ is contained inside the open ball $B_{\delta(\omega(0))}(\omega(0))$ with radius $\delta(\omega(0))$ and center $\omega(0)$.

\subsubsection{Manifold homotopially stratified spaces (MHSSs)}\label{S: MHSS}

A filtered space $X$ is a \emph{manifold homotopically stratified space (MHSS)} if the following conditions hold:
\begin{itemize}
\item $X$ is locally-compact, separable, and metric.

\item $X$ has finitely many strata, and each $X_i$ is an $i$-manifold without boundary and is locally-closed in $X$. 

\item For each $k>i$, $X_i$ is forward tame in $X_i\cup X_k$. 

\item For each $k>i$, the holink evaluation $\hl_s(X_i\cup X_k,X_i)\to X_i$ is a fibration. 

\item \label{I: CD} For each $x$, there is a stratum-preserving homotopy $\hl(X,x)\times I\to \hl(X,x)$ from the identity into a compact subset of $\hl(X,x)$.\footnote{This condition, requiring \emph{compactly dominated local holinks}, was not part of the original definition of Quinn \cite{Q1}. It first appears in the work of Hughes leading towards his Approximate Tubular Neighborhood Theorem in \cite{Hu02}. 
}  
\end{itemize}

We say that an MHSS $X$ is $n$-dimensional if its top manifold stratum has dimension $n$. This implies that $X$ is $n$-dimensional in the sense of covering dimension by \cite[Theorem III.2]{HuW}, which states that  a space that is the union of a countable number of closed subsets of dimension $\leq n$ has dimension $\leq n$. This condition holds for $X$ since each stratum is a separable manifold of dimension $\leq n$ (see also \cite[Theorem V.1]{HuW}). It then follows from \cite[Theorem III.1]{HuW} and \cite[Corollary II.16.34, Definition II.16.6, and Proposition II.16.15]{Br} that the cohomological dimension $\dim_{R} X$  of $X$ is $\leq n$ for any ring $R$ with unity (note that since $X$ is metric, it is paracompact). Similarly, $\dim_R Z\leq n$ for any subspace $Z$ of $X$. 

A subset of an MHSS is \emph{pure} if it is a closed union of components of strata. Each skeleto $X^i$ is a pure subset. The skeleton $X^{n-1}$ of an $n$-dimensional MHSS is also refered to as the \emph{singular set $\Sigma$}.

\subsection{Neighborhoods in stratified spaces}
\subsubsection{Teardrops}\label{S: teardrops}

 Given a map $p: X\to Y\times \R$, the \emph{teardrop} $X\cup_p Y$ of $p$ is the space $X\amalg Y$ with the minimal topology such that 
\begin{itemize}
\item $X\into X\cup_p Y$ is an open embedding, and
\item the function $c: X\cup_p Y\to Y\times (-\infty,\infty]$ defined by
\begin{equation*}
c(z)=\begin{cases}
p(z), &z\in X\\
(z,\infty), & z\in Y
\end{cases}
\end{equation*}
\end{itemize}
is continuous.

Given $f:X\to Y$, the teardrop $(X\times \R)\cup_{f\times \text{id}}Y$ is the \emph{open mapping cylinder of $f$ with the teardrop topology}. If $f$ is a proper map between locally compact Hausdorff spaces, then this is the usual mapping cylinder with the quotient topology (see \cite{Hu99a}). An alternative description of the teardrop topology of a mapping cylinder is as the topology  on $X\times (0,1)\amalg Y$  generated by the open subsets of $X\times (0,1)$ and sets of the form $U\cup (p^{-1}(U)\times (0, \epsilon))$, where $U$ is open in $Y$.

If $N$ is a nearly stratum-preserving deformation retract neighborhood (NSDRN) of a pure subset $Y$ of a manifold homotopically stratified space (MHSS), then $N$ is stratum-preserving homotopy equivalent to the mapping cylinder $M$ of the holink evaluation $\hl_s(N,Y)\to Y$, provided $M$ is given the teardrop topology. A proof can be found in \cite[Appendix]{GBF3}.

\subsubsection{Approximate tubular neighborhoods}

A \emph{weak stratified approximate fibration} $q:A\to B$ is a map such that the following lifting condition is satisfied: Given a diagram
\begin{diagram}
Z&\rTo^f & A\\
\dTo^{\times 0}&&\dTo_q\\
Z\times I&\rTo^F&B,
\end{diagram}
such that $Z$ is arbitrary and $F$ is a stratum-preserving homotopy, there is a \emph{weak stratified controlled solution} $\td F: Z\times I\times [0,1)\to A$ that is stratum-preserving along $I\times [0,1)$, satisfies $\td F(z,0,t)=f(z)$, and is such that the function $\bar F:Z\times I\times I\to B$ defined by $\bar F|Z\times I\times [0,1)=p\td F$ and $\bar F|Z\times I\times \{1\}=F$ is continuous.

A \emph{manifold stratified approximate fibration (MSAF)} is a proper map between MHSSs that is also a weak stratified approximate fibration. 
$N$ is  an approximate tubular neighborhood of a pure subset $Y$ of the MHSS $X$ if  there is a manifold stratified approximate fibration (MSAF) $p:N-Y\to Y\times \R$ such that the teardrop $(N-Y)\cup_p Y$ is homeomorphic to $N$.

The following Approximate Tubular Neighborhood Theorem is due to Hughes \cite{Hu02}, generalizing earlier special cases due to Hughes, Taylor, Weinberger, and Williams \cite{HTWW} and Hughes and Ranicki \cite{HR}:

\begin{theorem}[Approximate Tubular Neighborhood Theorem (Hughes)]
Let $X$ be a MHSS with compact singular set $\Sigma$ such that all non-minimal strata of $X$ have dimension $\geq 5$. If $Y\subset \Sigma$ is a pure subset of $X$, then $Y$ has an approximate tubular neighborhood in X. If $\Sigma$ is not compact, the theorem remains true if, in addition to the previous dimension restrictions, all noncompact strata are of dimension $\geq 5$ and the one-point compactification of $X$ is a MHSS with the point at infinity constituting a new stratum. 
\end{theorem}

By \cite[p. 873]{Hu02}, if $N$ is an approximate tubular neighborhood, then the natural extension of $p:N-Y\to Y\times \R$ to  $\td p: N\to Y\times (-\infty,\infty]$ is also an MSAF.

\section{Local approximate tubular neighborhoods}\label{S: LATN}

Let $X$ be a manifold homotopically stratified space (MHSS), and suppose $x$ is a point in the $k$th stratum $X_k=X^k-X^{k-1}$. We will say that $x$ has a \emph{local approximate tubular neighborhood} in $X$ if there is an open neighborhood $U$ of $x$ in $X$ such that 
\begin{enumerate}
\item $U\cap X^{k-1}=\emptyset$, 
\item $U\cap X_k\cong \R^k$, and 
\item $U$ is an approximate tubular neighborhood of $U\cap X_k$ in $(X-X^k)\cup (U\cap X_k)$. 
\end{enumerate}

We note that $(X-X^k)\cup (U\cap X_k)$ is an open subset of $X$ and thus is itself a MHSS according to \cite[Proposition 3.4]{GBF5}. Furthermore, $U\cap X_k$ is a closed union of components of strata in $(X-X^k)\cup (U\cap X_k)$ so that it is a pure subset. 

We say that the MHSS $X$ has \emph{sufficiently many local approximate tubular neighborhoods} if each point $x\in \Sigma$ possesses a local approximate tubular neighborhood. Any space for which Hughes's Approximate Tubular Neighborhood Theorem holds has  sufficiently many local approximate tubular neighborhoods:

\begin{proposition}\label{P: LATN}
Let $X$ be a manifold homotopically stratified space with compact singular space $\Sigma$ such that all the non-minimal strata of $X$ have dimension greater than or equal to five (or alternatively such that all the non-compact strata are of dimension greater than or equal to five and the one-point compactificiation of $X$ is a MHSS with the point at infinity constituting a new stratum).  Then $X$ has sufficiently many local approximate tubular neighborhoods.
\end{proposition}
\begin{proof}
By the Approximate Tubular Neighborhod Theorem \cite[Theorem 1.1]{Hu02} (see also \cite[Remark 7.2]{Hu02}), any pure subset $Y$ in $X$ has an approximate tubular neighborhood. We will use this to obtain local approximate tubular neighborhoods.

So suppose $x\in X_k$. Then $X^k$ itself is a pure subset, and we can suppose $X^k$ has an approximate tubular neighborhood $W$. So there is a manifold stratified approximate fibration (MSAF) $p:W-X^k\to X^k\times \R$ that extends continuously to $\bar p:W\to X^{k}\times (-\infty,\infty]$. But now let $V$ be a neighborhood of $x$ in $X_k$ homeomorphic to $\R^k$, and let $U=\td p^{-1}(V\times (-\infty,\infty])$. We claim that $U$ is a local approximate tubular neighborhood of $x$. 

It is apparent that $U$ is an open neighborhood of $x$ and that conditions (1) and (2) of the definition for a local approximate tubular neighborhood are satisfied. So we must check only that $U$ is an approximate tubular neighborhood of $V=U\cap X_k$. The restriction of $\td p$ to $U$ remains continuous, so we need only show that $p_U=p|U-V$ is an  MSAF. $U-V$ and $V\times \R$ are both MHSSs, and since $U-V=p^{-1}(V\times \R)$, $p_U$ is proper (the inverse image in $U-V$ of any compact set in $V\times \R$ is the same as its inverse image in $W-X^k$). Finally, we employ the fact that the restriction of any weak stratified approximate fibration to the inverse image of any open set is itself a weak stratified approximate fibration by the following lemma.
\end{proof}

\begin{lemma}\label{L: MSAF restrict}
Let $p:X\to Y$ be a weak stratified approximate fibration between metric spaces, and let $U$ be an open subset of $Y$. Then $p_U:p^{-1}(U)\to U$ is a weak stratified approximate fibration. 
\end{lemma}
\begin{proof}
We must show that $p_U$ possesses the weak stratified lifting property. So suppose we have a stratified lifting problem specified by $f:Z\to p^{-1}(U)$ and $F:Z\times I\to U$ for some metric space $Z$ (we may assume $Z$ to be metric by Remark 5.5b of \cite{Hug}). Certainly there is a weak stratified controlled solution $\td F:Z\times I\times [0,1)\to X$ such that $\td F(z,0,s)=f(z)$ for all $z,s$ and such that $\bar F:Z\times I\times I$ is continuous, where $\bar F=p\td F$ on $Z\times I\times [0,1)$ and $\bar F|_{Z\times I\times 1}=F$. We need to show that we can arrange for a new $\td F$ whose image is contained completely in $p^{-1}(U)$. By the continuity of $\bar F$, and since $\bar F(Z,I,1)\subset U$, there exists a neighborhood $W$ of $Z\times I\times 1$ in $Z\times I\times I$ such that $\bar F(W)\subset U$ and $\td F(W-Z\times I\times 1)\subset p^{-1}(U)$. Let $d:Z\times I\times1\to \R^{>0}$ be the distance from $(z,t,1)$ to $Z\times I\times I -W$. Now let $\td G:Z\times I\times [0,1) \to X$ be given by $\td G(z,t,s)\to \td F(z,t,s+ (1-s)(1-d(z,t)/2))$. $\td G$ maps to $p^{-1}(U)$ by construction, and it is a solution to the desired approximate lifting problem.
\end{proof}

\begin{lemma}\label{L: cofinal nghbds}
Suppose that $X$ is a MHSS and that $x\in X$ has a local approximate tubular neighborhood $U$. Then $x$ has a family $U=U_0\supset U_1\supset U_2\supset \cdots$ of local approximate tubular neighborhoods that is cofinal among all neighborhoods of $x$.
\end{lemma}
\begin{proof}
Suppose $x\in X_k$, and let $V=U\cap X_k$. Then $U$ is the teardrop of $p: U-V\to V\times \R$. Let $\td p:U\to V\times (-\infty,\infty]$ be the continuous extension. Since $U$ is a local approximate tubular neighborhood $V\cong \R^k$ by definition, and we may assume that $x$ is the origin of $\R^k$. Let $V_m=\frac{1}{m}D^k$, where $D^k$ is the open unit disk in $\R^k$. Let $W_m=V_m\times (m,\infty)$, and  let $U_m=\td p^{-1}(W_m\cup V_m)$. Then $U_m$ is certainly a neighborhood of $x$, and it follows from the same arguments  as used in the proof of Proposition \ref{P: LATN} that $U_m$ is an approximate tubular neighborhoods of $x$.    

To see that this system is cofinal, let $Z$ be any open neighborhood of $x$. We will show that some $U_m$ is a subset of $Z$. Suppose not. Then for all $m$, $U_m\cap (X-Z)\neq \emptyset$. So for each $W_M\cup V_m$, there is a point $x_m\in W_m\cup V_m$ such that $\td p^{-1}(x_m)\notin Z$. But we must have $\{x_m\}$ converge to $x$ and thus also $\td p^{-1}(x_m)$ converges to $x$, by definition of the teardrop topology. But then $x$ is a limit point of the closed set $X-Z$, a contradiction to $Z$ being an open neighborhoods of $x$. 
\end{proof}

\section{$\mc{IS^*}$ is the Deligne sheaf }\label{S: deligne}

In this section, we will demonstrate that if $X$ is an MHSS with sufficiently many approximate tubular neighborhoods, then the intersection chain sheaf $\mc{IS^*}$ is quasi-isomorphic to the Deligne sheaf $\mc P^*$. In \cite{GM2}, Goresky and MacPherson showed that the sheaf of simplicial intersection chains is quasi-isomorphic to $\mc P^*$ on PL pseudomanifolds. 
It was shown much later by the author in \cite{GBF10} that the singular intersection chain sheaf is quasi-isomorphic to the Deligne sheaf on topological pseudomanifolds. However,
it is by no means obvious that the Deligne sheaf and the intersection chain sheaf are quasi-isomorphic on MHSSs. On pseudomanifolds, one makes 
strong use of the geometric form of local distinguished neighborhoods $\R^{n-k}\times cL$, where $L$ is a compact pseudomanifold, but points in MHSSs have no such distinguished neighborhoods. It is true that each point has a neighborhood stratum-preserving homotopically equivalent to a space of the form $\R^{n-k}\times c\ms L$, and this fact is utilized in Quinn's proof of topological invariance of compactly supported intersection homology on MHSSs \cite{Q2}. But to establish the desired sheaf quasi-isomorphism, it is necessary  to consider closed support intersection homology on local neighborhoods, and these groups are not generally perserved under stratified homotopy equivalences (they would be preserved if the homotopy equivalences were proper, but in general they will not be - $\ms L$ generally will not even be locally compact). This necessitates the arguments to follow.

${}$\linebreak

Let $X$ be an $n$-dimensional MHSS with no codimension one stratum, let $\bar p$ be a fixed perversity, and let $\mc E$ be a local coefficient system on $X-X^{n-2}$. Let $U_k=X-X^k$, let $X_k=X^k-X^{k-1}$, and let $i_k:U_k\to U_{k+1}=U_k\cup X_{n-k}$ denote the inclusion.  We will omit $\bar p$  from the notation so long as it remains fixed. 

We recall that the Deligne sheaf $\mc P^*(\mc E)$ is defined inductively in  \cite{GM2} so that $\mc P^*_2=\mc E$ on $U_2=X-X^{n-2}$, and $$\mc P^*|_{U_{k+1}}=\mc P_{k+1}^*=\tau_{\leq \bar p(k)}Ri_{k*}\mc P_k^*$$ for $k\geq 2$. All formulas should be considered to live in the derived category of sheaves on $X$. In particular, $=$ really denotes quasi-isomorphism,  $Ri_{k*}$ is the derived functor of the sheaf pushforward $i_{k*}$, and $\tau_{\leq \bar p(k)}$ is the sheaf truncation functor. 

Let $\mc{IS}^*(\mc E)$ denote the sheaf of intersection chains on $X$ as defined in \cite{GBF10} with perversity $\bar p$ and local coefficients $\mc E$. We prove the following theorem:

\begin{theorem}\label{T: IS is deligne}
Let $X$ be an n-dimensional MHSS with no codimension one stratum  and with sufficiently many local approximate tubular neighborhoods. Let $\mc O$ be the orientation sheaf of the $n$-manifold $X-X^{n-2}$, and let $\mc E$ be a local coefficient system on $X-X^{n-2}$. Then the Deligne sheaf $\mc P^*(\mc E\otimes \mc O)$ and the sheaf of singular intersection chains $\mc{IS}^*(\mc E)$ are quasi-isomorphic. 
\end{theorem}

We begin by recalling the basic axioms of the Deligne sheaf in the version of \cite[Section V.2]{Bo}. Let $\ms S^*$ be a differential graded sheaf on the filtered space $X$ of finite cohomological dimension, and let $\ms S^*_k$ denote $\ms S^*|_{X-X^{n-k}}$. Then the axioms $AX1_{\bar p, X}\mc E$ are 
\begin{enumerate}
\item \label{I: basics} $\ms S^*$ is bounded, $\ms S^i=0$ for $i<0$, and $\ms S_2$ is the local coefficient system $\mc E$ on $X-X^{n-2}$,

\item \label{I: 0 range} For $k\geq 2$ and $x\in X_{n-k}$, $H^i(\ms S^*_x)=0$ if $i>\bar p(k)$, and 

\item \label{I: attaching range} The attachment map $\alpha_k:\ms S^*_{k+1}\into Ri_{k*}\ms S_k^*$ is a quasi-isomorphism up to (and including) $\bar p(k)$.
\end{enumerate}

We recall that for an inclusion $i:U\into V$, the attaching map of a sheaf complex $\ms S^*$ is simply the composition of the canonical sheaf maps $\ms S^*\to i_*i^*\ms S^*\to Ri_*i^*\ms S^*$. The first of these sheaf maps  corresponds to the restriction of sections of $\ms S^*$ to $U$, and the second map is induced by any quasi-isomorphism from $\ms S^*$ to a sheaf complex adapted to the functor $i_*$. It follows that on hypercohomology the attaching map induces a homomorphism $\H^*(V;\ms S^*)\to \H^*(U;\ms S^*)$ that can be interpreted as being induced by the restriction of sections of any appropriate resolution of $\ms S^*$. If $x\in X_{n-k}$, Axiom \eqref{I: attaching range} is therefore equivalent to the condition that $H^i(\ms S^*)_x\cong \varinjlim_{x\in U}\H^i(U-U\cap X_{n-k};\ms S^*)$ for all $i\leq \bar p(k)$ (see \cite[V.1.7]{Bo}, \cite[Section 4.1.4]{BaIH}).

By \cite[Theorem V.2.5]{Bo}, any sheaf that satisfies the axioms $AX1_{\bar p,X}\mc E$ is quasi-isomorphic to $\mc P^*(\mc E)$. This theorem is stated for pseudomanifolds, but the proof applies for any filtered space. Thus we set out to show that $\mc{IS}^*$ satisfies the axioms.

As observed in \cite[Section V.2.7]{Bo}, since we are really working in the derived category, the first two conditions of axiom \eqref{I: basics} can be replaced with the conditions that $\ms S^*$ is bounded below and that $\mc H^i(\ms S^*)=0$ for $i<0$ and for $i\gg 0$.
And, in fact, the strict bounded below condition is never used in the proof of \cite[Section V.2.7]{Bo}; it seems to be invoked only later in \cite[Remark V.2.7.b]{Bo} to assure the convergence of the hypercohomology spectral sequence. Since we noted in \cite{GBF10} that the hypercohomology spectral sequence of $\mc{IS}^*$ does indeed converge with no difficulty (since $\mc{IS^*}$ is homotopically fine), there is both no such difficulty here \emph{and} this condition is unnecessary to prove the desired quasi-isomorphism. Thus it suffices to demonstrate that $\mc{IS}^*$ satisfies the axioms, except for the strict bounded below condition. (Additionally, once we have shown that $\mc{IS}^*$ satisifies the other properties, we can note that the condition $\mc H^i(\mc{IS}^*)=0$ for $i<0$ implies that $\mc{IS}^*$ is quasi-isomorphic to $\tau_{\geq 0}\mc{IS}^*$, which then itself satisfies all of the axioms).

We begin with axiom \eqref{I: 0 range}. 

\begin{proposition}\label{P: 0 range}
Let $\mc{IS}^*$ be the intersection chain sheaf on the MHSS $X$ with sufficiently many local approximate tubular neighborhoods. Then for $k\geq 2$ and $x\in X_{n-k}$, $H^i(\mc{IS}^*_x)=0$ if $i>\bar p(k)$ (i.e. $\mc{IS^*}$ satisfies Axiom \eqref{I: 0 range}). 
\end{proposition} 
\begin{proof}
Suppose $x\in X_{n-k}$. By elementary sheaf theory, $H^i(\mc{IS}^*_x)=\varinjlim_{x\in U}\H^i(U;\mc{IS}^*)$. By \cite[Proposition 3.7]{GBF10}, the restriction of $\mc{IS}^*$ to an open subset is quasi-isomorphic to the intersection chain sheaf on the subset, and thus $\H^i(U;\mc{IS}^*)\cong IH_{n-i}^{\infty}(U)$. 

Now, suppose that $U$ is a local approximate tubular neighborhood of $x$ (and hence an approximate tubular neighborhood of $U\cap X_{n-k}\cong \R^{n-k}$). By \cite[Corollary 9.2 and Proposition 9.34]{GBF13}, $U$ is also an outwardly stratified tame nearly stratum-preserving deformation retract neighborhood of $U\cap X_{n-k}\cong \R^{n-k}$ (the reader who wants to know what all that means is urged to consult \cite{GBF13}; we will merely use this fact to invoke some other results from \cite{GBF13} regarding such neighborhoods). 

Putting together Theorem 6.15 and Proposition 5.1 of \cite{GBF13}, since $U$ is an outwardly stratified tame nearly stratum-preserving deformation retract neighborhood, $IH^{\infty}_{n-*}(U)$ is the abutment of a spectral sequence with $E_2$ terms $E_2^{r,s}\cong H^r(\R^{n-k};IH^c_{n-(n-k)-s}(c\ms L,\ms L\times \R))$, where here $\ms L=\hl_s(U,x)$. The coefficient system is constant because the base space is homeomorphic to $\R^{n-k}$. We note also that the pair $(c\ms L,\ms L\times \R)$ is the pre-image of the point $x$ under the collapse map $(M,M-\R^{n-k})\to \R^{n-k}$, where $M$ is the mapping cylinder of the holink evaluation $\pi: \hl_s(U,\R^{n-k})$. This mapping cylinder, and hence also the cone, are given the teardrop topology. 

This spectral sequence collapses immediately, all terms being $0$ except for the terms $E_2^{0,s}$, at which we have $E_2^{0,s}\cong IH^c_{k-s}(c\ms L,\ms L\times \R)$. So we have $\H^i(U;\mc{IS}^*)\cong IH_{n-i}^{\infty}(U)\cong E_2^{0,i}\cong IH^c_{k-i}(c\ms L,\ms L\times \R)$. Now, $\ms L$ is an infinite dimensional space, but under the conventions for interesection homology under stratified homotopy equivalences (see, e.g, \cite{GBF13}), all strata of $M$ and $\ms L\times \R\subset c\ms L\subset M$ simply inherit the formal dimension labels of the strata they arise from in $U$ under the stratum-preserving homotopy equivalence $M\sim_{sphe}U$ (see Section \ref{S: teardrops}). At the same time, it is not these formal dimensions that really matter in intersection homology, only the codimensions, so we are free to shift all dimension labels. If we subtract $n-k$ from all the formal strata dimensions on $(c\ms L, \ms L\times \R)$, then the cone point has dimension $0$, as appropriate, and we see that we are free to apply the usual intersection homology cone formula (see \cite[Proposition 5]{Ki}), as the standard arguments of its proof will apply with $c\ms L$ a filtered space of filtered dimension $k$.  Thus 
\begin{equation*}
IH^c_j(c\ms L,\ms L\times \R)\cong 
\begin{cases}
0, & j\leq k-1 -\bar p(k),\\
IH^c_{j-1}(\ms L), &j>k-\bar p(k).
\end{cases}
\end{equation*}
Thus $\H^i(U;\mc{IS}^*)=0$ for $k-i\leq k-1-\bar p(k)$, i.e. for $i>\bar p(k)$.

Since $x$ possesses a cofinal system of local approximate tubular neighborhoods, the proposition follows.

\end{proof}

Next we start work towards the proof that $\mc{IS^*}$ satisfies axiom \eqref{I: attaching range}.

\begin{lemma}\label{L: out tame}
Let $p: X\to \R^m$ (or $p:X\to \R^{m}_+$, the closed half space) be a proper weak stratified approximate fibration. Then $IH^{\infty}_*((X-p^{-1}(0))_X)=0$, and hence $IH^{\infty}_*(X)\cong IH^{\infty}_*(X, X-p^{-1}(0))$.  
\end{lemma}
\begin{proof}
The second statement follows from the long exact sequence of the pair, once we prove the first. The proof that $IH^{\infty}_*((X-p^{-1}(0))_X)$ comes by showing that we can ``push cycles off to infinity''. A very similar statement and proof can be found in Proposition 6.7 of \cite{GBF13}, though in  a slightly different context.  The main point is that we need to show that 
$X-p^{-1}(0)$ possesses a version of the property that we refer to in \cite{GBF13} as ``outward stratified tameness'' of approximate tubular neighborhoods. The definition of this term in \cite[Section 6.2]{GBF13} applies to certain neighborhoods, but the appropriate modified condition here would say that for any metric space $Z$ and any proper map $g:Z\to X$ such that $g(Z)\in X-p^{-1}(0)$, there exists a propert stratum-preserving homotopy $H:Z\times [0,\infty)\to X$ such that $H(Z\times [0,\infty))\subset X-p^{-1}(0)$ and $H|_{Z\times 0}=f$. Once this condition is established, the proof that the intersection homology groups are $0$ follows by a direct modification of the proof of the cited proposition. We let the reader consult that proof for precise details; the idea is  
that for any intersection cycle $\xi$, this outward tamesness property allows us to build a homotopy of $|\xi|$ out to infinity. Then, this homotopy is used to build the desired infinite-chain null-homology.

Thus we should concentrate on the proof of existence of such proper open-ended homotopies $H$. Here, also, the proof is very similar to that of the proof that approximate tubular neighborhoods are outwardly stratified tame \cite[Proposition 9.3]{GBF13}. The proof of that proposition is rather lengthy, so we will not reproduce a modifed version here. The interested reader should note that the appropriate modifcation is to replace the sets $K_i\times (-i,\infty)$ with the closed disks (or half-disks for $\R^{m}_+$) $D_i$ of radius $i$ in $\R^m$, and the sets $C_i$ with $p^{-1}(D_i)$. Then one proceeds as in that proof to construct $H$ so that $H(\cdot,0)=f$ and for each positive integer $i$,  $H(Z\times [i,\infty))\subset X-C_i$, and $H(z,t)=f(z)$ if $t\in[0,i]$ and $z\in f^{-1}(X-\text{int}(C_{i+1}))$ (this last condition is the key to properness, since at each finite time only a compact set is moved by the homotopy). The proofs that such an $H$ suffices and that it can be constructed are similar to those of \cite[Proposition 9.3]{GBF13}, and we leave the necessary modifications to the reader. 
\end{proof}

\begin{corollary}\label{C: del ngbd reduce}
Let $U$ be a local approximate tubular neighborhood of the point $x$ in the stratum $X_{n-k}$ of the MHSS $X$. Let $V=U-U\cap X_{n-k}$, and let  $p:V\to (U\cap X^{n-k})\times \R$ be the proper MSAF of the definition of approximate tubular neighborhoods. Let $y$ be any point in $(U\cap X^{n-k})\times \R$. Then $IH^{\infty}_*(V)\cong IH^{\infty}_*(V, V-p^{-1}(y))\cong IH^c_*(V,V-p^{-1}(y))\cong IH^c_*(U,U-p^{-1}(y))$.  
\end{corollary}
\begin{proof}
By the definition of local approximate tubular neighborhoods, $U\cap X^{n-k}\cong \R^{n-k}$. Thus $(U\cap X^{n-k})\times \R\cong \R^{n-k+1}$. We can  treat  any $y$ as the origin of $\R^{n-k-1}$ and apply the preceding lemma to obtain the first isomorphism. Since $y$ is compact and $p$  is proper, $V-p^{-1}(y)$ is cocompact, and so  the second isomorphism follows by \cite[Lemma 2.12]{GBF10}. The last isomorphism is by excision (see \cite[Lemma 2.11]{GBF10}).  
\end{proof}

\begin{corollary}\label{C: ngbd reduce}
Let $U$ be a local approximate tubular neighborhood of the point $x$ in the stratum $X_{n-k}$ of the MHSS $X$. Let  $\td p:U\to (U\cap X^{n-k})\times (-\infty,\infty]$ be the proper MSAF arising from the definition of approximate tubular neighborhoods. Let $(y,t)$ be any point in $(U\cap X^{n-k})\times (-\infty,\infty)$. Then $IH^{\infty}_*(U)\cong IH^{\infty}_*(U, U-\td p^{-1}(y\times [t,\infty]))\cong IH^c_*(U, U-\td p^{-1}(y\times [t,\infty]))$.  
\end{corollary}
\begin{proof}
By the definition of local approximate tubular neighborhoods, $U\cap X^{n-k}\cong \R^{n-k}$. Thus $(U\cap X^{n-k})\times (-\infty,\infty]\cong \R^{n-k+1}_+$, and we can treat $y\times\infty$ as the origin in $\R^{n-k+1}_+$. Note also that $\td p^{-1}(y\times \infty)$ is just a single point in $X_{n-k}$, which we can also call $y\times \infty$. Thus by Lemma \ref{L: out tame},  $IH^{\infty}_*(U)\cong IH^{\infty}_*(U, U-(y\times \infty])) \cong IH^{\infty}_*(U, U-\td p^{-1}(y\times \infty]))$. By \cite[Lemma 2.12]{GBF10}, this is isomorphic to $IH^{c}_*(U, U-(y\times \infty])) \cong IH^c_*(U, U-\td p^{-1}(y\times \infty]))$. Finally, we see that  this is isomorphic to $IH^c_*(U, U-\td p^{-1}(y\times [t,\infty]))$, by the long exact sequence of the triple, since $IH_*^c(U-\td p^{-1}(y\times [t,\infty]), U-\td p^{-1}(y\times \infty))=0$: clearly $\R^{n-k+1}_+-y\times \infty$ deformation retracts into $\R^{n-k+1}_+-y\times [t,\infty]$, and this may be used to push chains around appropriately in $U$, using the MSAF $\td p$. 
\end{proof}

\begin{proposition}
Let $U$ be a local approximate tubular neighborhood of the point $x$ in the stratum $X_{n-k}$ of the MHSS $X$. Let  $\td p:U\to (U\cap X^{n-k})\times (-\infty,\infty]$ be the proper MSAF arising from the definition of approximate tubular neighborhoods. Then $IH^c_*(U, U-\td p^{-1}(x\times [t,\infty]))\to IH^c_*(U,U-\td p^{-1}(x\times t))$, induced by inclusion, is an isomorphism for $*\geq n-\bar p(k)$. 
\end{proposition}
\begin{proof}
To simplify the notation, recall that  $U\cap X_{n-k}\cong \R^{n-k}$ and assume, without loss of generality, that $x=0\in \R^{n-k}$. We will use the long exact sequence of the triple $(U,  U-\td p^{-1}(0\times t), U-\td p^{-1}(0\times [t,\infty]))$. 

So consider $IH_*^c( U-\td p^{-1}(0\times t), U-\td p^{-1}(0\times [t,\infty]))$. By excision, this is isomorphic to $IH_*^c(\td p^{-1}(\R^{n-k}\times (t,\infty]), \td p^{-1}((\R^{n-k}-0)\times (t,\infty]))$. Using Lemma \ref{L: MSAF restrict}, the restriction of a proper MSAF to the inverse image of an open subset is again a proper MSAF, so   $\td p^{-1}(\R^{n-k}\times (t,\infty])$ and $\td p^{-1}((\R^{n-k}-0)\times (t,\infty])$ are approximate tubular neighborhoods respectively of $\R^{n-k}$ and $\R^{n-k}\times 0$. In particular, then, by \cite[Section 9]{GBF13}, these neighborhoods are each stratum-preserving homotopy equivalent to mapping cylinders of homotopy link evaluations.

Let $M$ be the mapping cylinder of the holink evaluation $\hl_s(\td p^{-1}(\R^{n-k}\times (t,\infty]),\R^{n-k})\to \R^{n-k}$, which is a stratified fibration by \cite{Hug}, and let $P$ be the mapping cylinder collapse, which is also a stratified fibration by \cite[Proposition 3.3]{GBF5}. 
We will show below that $IH_*^c(\td p^{-1}(\R^{n-k}\times (t,\infty]), \td p^{-1}((\R^{n-k}-0)\times (t,\infty]))\cong IH^c_*(M, M-P^{-1}(0))$. Let us assume for now that this isomorphism holds. Then, using  \cite[Corollary 3.4]{GBF3}, there is  a stratum- and fiber-preserving homotopy equivalence from $M$ to $\R^{n-k}\times F$, where $F=P^{-1}(0)$. Since $M$ is a mapping cylinder of $\pi$, $F=P^{-1}(0)=c\pi^{-1}(0)=c\hl_s(\td p^{-1}(\R^{n-k}\times (t,\infty]),0)$, where $c$ indicates the open cone. Thus $IH^c_*(M, M-P^{-1}(0))\cong IH^c_*((\R^{n-k}, \R^{n-k}-0)\times F)$, which, employing the intersection homology K\"unneth theorem (which is allowed since $(\R^{n-k}, \R^{n-k}-0)$ is a manifold pair), is homeomorphic to $IH_{*-(n-k)}^c(F)$. Since $F$ is the cone on $\hl_s(\td p^{-1}(\R^{n-k}\times (t,\infty]),0)$, we may argue again as in Proposition \ref{P: 0 range} to conclude that we can employ the standard cone formula as though $\hl_s(\td p^{-1}(\R^{n-k}\times (t,\infty]),0)$ were a $k-1$ dimensional space. 

Thus 
\begin{equation*}
IH^c_j(F)\cong
\begin{cases}
0,& j\geq k-1-\bar p(k),\\
IH^c_j(\hl_s(\td p^{-1}(\R^{n-k}\times (t,\infty]),0)), & j<k-1-\bar p(k).
\end{cases}
\end{equation*}
So 
$IH_*^c( U-\td p^{-1}(0\times t), U-\td p^{-1}(0\times [t,\infty]))\cong IH_{*-(n-k)}^c(F)=0$ if $*\geq n-1-\bar p(k)$, and by the long exact sequence of the triple,
$IH^c_*(U, U-\td p^{-1}(x\times [t,\infty]))\cong IH^c_*(U,U-\td p^{-1}(x\times t))$ for $*\geq n-\bar p(k)$, as desired.

It remains to show that  $IH_*^c(\td p^{-1}(\R^{n-k}\times (t,\infty]), \td p^{-1}((\R^{n-k}-0)\times (t,\infty]))\cong IH^c_*(M, M-P^{-1}(0))$. The proof is similar to some of those in \cite{GBF13}: Let $\delta:\R^{n-k}-0\to (0,\infty)$ be a continuous function such that for $z\in \R^{n-k}-0$, $\delta(z)$ is less than the distance in $X$ from $z$ to $X-\td p^{-1}((\R^{n-k}-0)\times (t,\infty])$. Let $M^{\delta}_0$ be the mapping cylinder of the evaluation $\hl_s^{\delta}(\td p^{-1}((\R^{n-k}-0)\times (t,\infty]),\R^{n-k}-0)\to \R^{n-k}-0$. Then the inclusion $(M, M_0^\delta)\into (M, M-P^{-1}(0))$ induces a stratum-preserving homotopy equivalence $M_0^{\delta}\to M-P^{-1}(0)$ by the arguments of Quinn \cite{Q1}. So, employing the five lemma, $IH_*^c(M, M_0^\delta)\cong IH^c_*(M, M-P^{-1}(0))$.

On the other hand,  by \cite{GBF13}, the approximate tubular neighborhood  $\td p^{-1}(\R^{n-k}\times (t,\infty])$  is a nearly stratum-preserving deformation retract neighborhood, and so by  \cite[Proposition A.1]{GBF3}, it is stratum-preserving homotopy equivalent $M$. If we let $g$ be the modified path evaluation map of the proof of \cite[Proposition A.1]{GBF3}, we see that $g$ maps the pair $(M, M_0^{\delta})$ to 
the pair $(\td p^{-1}(\R^{n-k}\times (t,\infty]), \td p^{-1}((\R^{n-k}-0)\times (t,\infty]))$. But $g:M\to \td p^{-1}(\R^{n-k}\times (t,\infty])$ is precisely the stratum-preserving homotopy equivalence of the cited proposition.  The restriction $g: M_0^{\delta} \to \td p^{-1}((\R^{n-k}-0)\times (t,\infty]) $ is also a stratum-preserving homotopy equivalence since it factors as the composition of two stratum-preserving homotopy equivalences $M_0^{\delta}\to M_0 \to \td p^{-1}((\R^{n-k}-0)\times (t,\infty])$, where $M_0$ is the mapping cylinder of the holink evaluation $\hl_s( \td p^{-1}((\R^{n-k}-0)\times (t,\infty]), \R^{n-k}-0)\to \R^{n-k}-0$, the first map is inclusion, which is a stratum-preserving homotopy equivalence by Quinn \cite{Q1}, and the second map is again the homotopy equivalence of \cite[Proposition A.1]{GBF3}. Note that the claimed factorization holds, since we may choose compatible shrinking maps $S$, as defined in the proof of \cite[Proposition A.1]{GBF3}, for all involved holink spaces. So,  employing the five-lemma, $IH_*^c(M, M_0^\delta)\cong IH_*^c(\td p^{-1}(\R^{n-k}\times (t,\infty]), \td p^{-1}((\R^{n-k}-0)\times (t,\infty]))$. This completes the proof. 

\begin{comment}
That $IH_*^c(\td p^{-1}(\R^{n-k}\times (t,\infty]))\cong IH^c_*(M)$ is immediate, since we have noted that $\td p^{-1}(\R^{n-k}\times (t,\infty])$ is an approximate tubular neighborhood of $\R^{n-k}$ and thus by \cite[Corollary 9.2]{GBF13}, it is a nearly startum-preserving deformation retract neighborhood or $\R^{n-k}$, and by  \cite[Proposition A.1]{GBF3} (which corrects a proof of Quinn in \cite{Q1}, $\td p^{-1}(\R^{n-k}\times (t,\infty])$ is stratum-preserving homotopy equivalent to the mapping cylinder of the holink evaluation $\pi: \hl_s(\td p^{-1}(\R^{n-k}\times (t,\infty]),\R^{n-k})\to \R^{n-k})$. Let $f: \td p^{-1}(\R^{n-k}\times (t,\infty])\to M$ be the stratum-preserving homotopy equivalence. To complete the proof, it suffices by the 
\end{comment}
\end{proof}

\begin{corollary}\label{C: ngb iso}
Let $U$ be a local approximate tubular neighborhood of the point $x$ in the stratum $X_{n-k}$ of the MHSS $X$. Then the restriction $IH_*^{\infty}(U)\to IH_*^{\infty}(U-U\cap X_{n-k})$ is an isomorphism for $*\geq n-\bar p(k)$. 
\end{corollary}
\begin{proof}
Applying the preceding proposition, it suffices to show that the following diagram commutes:
\begin{equation*}
\begin{CD}
IH_*^{\infty}(U)&@>>>& IH_*^{\infty}(U-U\cap X_{n-k})\\
@V\cong VV&&@VV\cong V\\
IH^c_*(U, U-\td p^{-1}(x\times [t,\infty]))&@>>>& IH^c_*(U,U-\td p^{-1}(x\times t))
\end{CD}
\end{equation*}
The vertical maps are the isomorphisms of Corollaries  \ref{C: del ngbd reduce} and \ref{C: ngbd reduce}. But the commutativity is easy to see at the chain level using representative cycles and the techniques of, for example, the proof of \cite[Lemma 2.12]{GBF10}: One begins with a chain $\xi$ representing a cycle in $IC^{\infty}_*(U)$ and ends up with a relative cycle in $IC_*^c(U,U-\td p^{-1}(x\times t))$ that is obtained by sufficiently subdividing $\xi$ and then excising all but a finite number of singular simplices whose supports lie in $U-\td p^{-1}(x\times t)$. 
We are not free to perform excisions on intersection chains along just any boundaries of simplices but a procedure for performing allowable excisions of intersection chains was well-established in \cite{GBF10} and may be applied here. By considering what happens to $\xi$, one sees that the two different ways of chasing around the diagram  yield the same result.  
\end{proof}

\begin{proposition}\label{P: attaching range}
Let $\mc{IS}^*$ be the intersection chain sheaf on the MHSS $X$ with sufficiently many local approximate tubular neighborhoods. Then $\alpha_k:\mc{IS}^*_{k+1}\into Ri_{k*}\mc{IS}_k^*$ is a quasi-isomorphism up to (and including) $\bar p(k)$ (i.e. $\mc{IS^*}$ satisfies Axiom \eqref{I: attaching range}). 
\end{proposition} 
\begin{proof}
Let $x\in X_{n-k}$. Then $H^i(\mc{IS}^*_x)=\varinjlim_{x\in U}\H^i(U;\mc{IS}^*)\cong \varinjlim_{x\in U}IH^{\infty}_{n-i}(U)$. For the last isomorphism, we use that  the restriction of $\mc{IS}^*$ to an open subset is quasi-isomorphic to the intersection chain sheaf on the subset  by \cite[Proposition 3.7]{GBF10}. Similarly, $H^i(Ri_{k*}\mc{IS}_k^*)_x=\varinjlim_{x\in U}\H^i(U; Ri_{k*}\mc{IS}_k^*)\cong \varinjlim_{x\in U}\H^i(U-U\cap X_{n-k}; \mc{IS}^*)\cong \varinjlim_{x\in U} IH_{n-i}^{\infty}(U-U\cap X_{n-k})$. By Corollary \ref{C: ngb iso}, for a fixed $U$, $IH^{\infty}_{n-i}(U)\cong H_{n-i}^{\infty}(U-U\cap X_{n-k})$ for $i\leq \bar p(k)$, where the isomorphism is induced by restriction. But this restriction is compatible with the attaching map (see the discussion of the attaching map following the statement of Theorem \ref{T: IS is deligne}), and it suffices to show that if $V\subset U$ is another approximate tubular neighborhood of $x$ from the cofinal system of Lemma \ref{L: cofinal nghbds}
then the vertical maps are isomorphisms in the following commutative diagram, in which all maps are induced by restrictions:
\begin{equation*}
\begin{CD}
IH^{\infty}_{n-i}(U) &@>>>& H_{n-i}^{\infty}(U-U\cap X_{n-k})\\
@V\cong VV&&@VV\cong V\\
IH^{\infty}_{n-i}(V) &@>>>& H_{n-i}^{\infty}(V-V\cap X_{n-k}).
\end{CD}
\end{equation*}

For the lefthand vertical isomorphism, let us identify $U\cap X_{n-k}$ with  $\R^{n-k}$, let $x=0$, and let $\td p: U-\R^{n-k} \to \R^{n-k}\times (-\infty, \infty]$ be the proper MSAF of the definition of the approximate tubular neighborhood. We may suppose as in Lemma \ref{L: cofinal nghbds} that $V\cong \td p^{-1}(D_m\times [m, \infty])$, where $D_m$ is the disk of radius $m$ in $\R^{n-k}$. Let $y=0\times (m+1)\in \R^{n-k}\times (-\infty,\infty]$. Then $IH^{\infty}_{n-i}(U)\cong IH^c_{n-i}(U,U-p^{-1}(y))$ by Corollary \ref{C: del ngbd reduce}, and similarly $IH^{\infty}_{n-i}(V)\cong IH^c_{n-i}(V,V-p^{-1}(y))$. But $IH^c_{n-i}(U,U-p^{-1}(y))\cong IH^c_{n-i}(V,V-p^{-1}(y))$ by excision, and once again by using small enough representative cycles (which we can choose to have support arbitrarily close to $p^{-1}(y)$ by using sufficiently fine subdivisions and excisions (see \cite{GBF10})), this isomorphism is compatible with restriction of infinite chains. The proof for the righthand vertical maps follows similarly from   Corollary \ref{C: ngbd reduce}.
\end{proof}

Finally, we attack Axiom \ref{I: basics}.

\begin{proposition}\label{P: top stratum}
Let $X$ be an MHSS, and let $\mc{IS}^*$ be the intersection chain sheaf with coefficients in the local system $\mc E$ on $X-X^{n-2}$. Then $\mc{IS}^*|_{X-X^{n-2}}$ is quasi-isomorphic to $\mc E\otimes \mc O$.
\end{proposition}
\begin{proof}
By Proposition 3.7 of \cite{GBF10}, the restriction of $\mc{IS^*}$ to $X-X^{n-2}$ is quasi-isomorphic to the intersection chain sheaf on $X-X^{n-2}$. But $X-X^{n-2}$ is a manifold, so $\mc{IS}^*$ is simply the ordinary singular chain sheaf, whose local cohomology groups are $H^i(\mc{IS}^*_x)\cong H_{n-i}(X, X-x;\mc E)$. The proposition follows. 
\end{proof}

\begin{proposition}\label{P: basics}
Let $X$ be an MHSS, and let $\mc{IS}^*$ be the intersection chain sheaf on $X$. Then $H^i(\mc{IS}^*_x)=0$ for $i<0$ and $i>n$. 
\end{proposition}
\begin{proof}
Since $\mc{IS}^*$ is the sheafification of the presheaf $U\to IC_{n-*}(X,X-\bar U)$, it is immediate that $\mc{IS}^i$ is identically $0$ for $i>n$. 

For $i<0$, we will induct down over the strata of $X$, starting with the top stratum $X-X^{n-2}$. Since this stratum is a manifold and open in $X$, the restriction to it of $\mc{IS}^*$ is quasi-isomorphic to the ordinary singular chain sheaf, using  \cite[Propisition 3.7]{GBF10}. Thus $H^i(\mc{IS}^*_x)\cong H^{\infty}_{n-i}(U)$, where $U$ is any euclidean neighborhood of $x$. This is certainly $0$ for $i<0$. 

Now, we assume by induction that $H^i(\mc{IS}^*_z)=0$ for $i<0$ if $z\in X-X^{n-k}$ and show that $H^i(\mc{IS}^*_x)=0$ for $i<0$ if $x\in X_{n-k}=X^{n-k}-X^{n-k-1}$. By Proposition \ref{P: attaching range}, $\mc H^i(\mc{IS}^*_{k+1})_x\cong \mc H^i( Ri_{k*}\mc{IS}_k^*)_x$ for  $i\leq \bar p(k)$, and since $\bar p(k)\geq 0$, this applies for $i<0$. We have  $\mc H^i( Ri_{k*}\mc{IS}_k^*)_x\cong \varinjlim_{x\in U}\H^i(U;Ri_{k*}\mc{IS}_k^*)\cong 
\varinjlim_{x\in U}\H^i(U-U\cap X_{n-k};\mc{IS}_k^*)$.  But now $\H^i(U-U\cap X_{n-k};\mc{IS}_k^*)$ is the abutment of the hypercohomology spectral sequence with $E_2$ terms $E_2^{p,q}\cong H^p(U-U\cap X_{n-k}; \mc H^q(\mc{IS}_k^*))$, and by the induction hypothesis, these groups are $0$ if either $q$ or $p$ is $<0$. So then any term of the spectral sequence that would contribute to $\H^i(U-U\cap X_{n-k};\mc{IS}_k^*)$ for $i<0$ is trivial, and all these hypercohomology groups are $0$. Thus $H^i(\mc{IS}^*_x)=0$ for $i<0$ if $x\in X_{n-k}$, and the proof is completed by induction.
\end{proof}

As noted in our discussion of the axioms following the statement of Theorem \ref{T: IS is deligne}, Propositions \ref{P: 0 range}, \ref{P: attaching range}, \ref{P: top stratum} and \ref{P: basics} suffice to demonstrate that $\mc{IS}^*$ and the Deligne sheaf are quasi-isomorphic, proving the theorem. We note also that we can make $\mc{IS}^*$ legitimately bounded below by replacing it with $\tau_{\geq m}\mc{IS}^*$ for any $m\leq 0$: It follows from Proposition \ref{P: basics} that $\mc{IS}^*$ and $\tau_{\geq m}\mc{IS}^*$ are quasi-isomorphic and so certainly $\tau_{\geq m}\mc{IS}^*$ satisfies Axioms \eqref{I: basics} and \eqref{I: 0 range}. But also $\mc H^i(Ri_{k*}\tau_{\geq m}\mc{IS}_k^*)_x\cong\varinjlim_{x\in U}\H^i(U; Ri_{k*}\tau_{\geq m}\mc{IS}_k^*)\cong \varinjlim_{x\in U}\H^i(U-U\cap X_{n-k}; \tau_{\geq m}\mc{IS}^*)\cong \varinjlim_{x\in U}\H^i(U-U\cap X_{n-k};\mc{IS}^*)\cong \mc H^i(Ri_{k*}\mc{IS}_k^*)_x$, the next to last isomorphism since $\tau_{\geq m}\mc{IS}^*$ and $\mc{IS}^*$ are quasi-isomorphic. Thus, employing Proposition \ref{P: attaching range} and once again the quasi-isomorphism of $\tau_{\geq m}\mc{IS}^*$ and $\mc{IS}^*$, $\mc H^i(Ri_{k*}\tau_{\geq m}\mc{IS}_k^*)_x$ and $\tau_{\geq m}\mc{IS}^*_{k+1}$ are quasi-isomorphic in the appropriate range. Thus $\tau_{\geq m}\mc{IS}^*$ satisfies \emph{all} of the axioms spot on, including the boundedness, and is quasi-isomorphic to the Deligne sheaf. Thus if desired (though not necessary), we can use $\tau_{\geq m}\mc{IS}^*$ as a bounded below intermediary that is quasi-isomorphic  both to the Deligne sheaf \emph{and} to $\mc{IS}^*$. \hfill\qedsymbol

\subsection{Superperversities and codimension one strata}

For simplicity of the preceding discussion, we have assumed the MHSS $X$ to have no codimension one statum, i.e. $X^{n-1}-X^{n-2}=\emptyset$, and we have assumed the perversity $\bar p$ to be traditional, i.e. $\bar p(1)=\bar p(2)=0$. These restrictions comply with  those originally imposed by Goresky and MacPherson in their initial development of intersection homology theory \cite{GM1, GM2}. However,
it became apparent in the work of Cappell and Shaneson, particularly in their Superduality Theorem \cite{CS}, that it is also fruitful in the context of Deligne sheaf intersection homology to study \emph{superperversities} - those perversities $\bar p$  for which $\bar p(1)$ or $\bar p(2)$ is greater than $0$ (though we still require that $\bar p(k)\leq \bar p(k+1)\leq \bar p(k)+1$).\footnote{By contrast, any consideration of Deligne sheaves using \emph{subperversities} trivializes immediately, since truncation $\tau_{\leq \bar p(k)}$ yields the $0$ sheaf complex if $\bar p(k)<0$.} Superperverse intersection homology has since been studied in a variety of other contexts; see, e.g., \cite{HS91, GBF10, GBF11, Sa05}.

It was shown in \cite{GBF10} that if $X$ is a pseudomanifold, possibly with codimension one stratum, and $\bar p$ is a superperversity, then the Deligne sheaf intersection homology is isomorphic to singular chain intersection homology with coefficients in a certain \emph{stratified coefficient system} $\mc E_0$ based on the coefficient system $\mc E$. By definition, the simplices in these singular chains carry coefficients over the subsets of their supports that do not intersect $X^{n-1}$, and they carry a formal $0$ coefficient over the subsets of their supports that do intersect $X^{n-1}$. The reader is advised to consult \cite{GBF10} for further details, but the point is that this coefficient convention allows for $1$-chains whose boundary $0$-chains lie in $X^{n-1}$. This manages to correct a technical deficiency that otherwise prevents sheaf-theoretic and singular chain-theoretic intersection homology from agreeing. In particular, with this coefficient correction, superperverse singular chain intersection homology satisfies the usual intersection homology cone formula, and the  superperverse singular intersection chain sheaf satisifies the Deligne sheaf axioms.

Returning now to MHSSs,  our proof of Theorem \ref{T: IS is deligne} holds even if $X$ has a codimension one stratum or $\bar p$ is a superperversity, provided we replace $\mc{IS}^*$ with the sheaf of singular intersection chains with stratified coefficients as in \cite{GBF10}. All of the arguments we have employed involving 
excision, subdivision, K\"unneth theorems, and stratum-preserving homotopy invariance of compactly supported intersection homology hold for this variant (see \cite{GBF10}), as well as the cone formulas, which are at the crux of all the computational arguments. Thus we can generalize Theorem \ref{T: IS is deligne} as follows:

\begin{theorem}
Let $X$ be an n-dimensional MHSS, possibly  with codimension one stratum  and with sufficiently many local approximate tubular neighborhoods. Let $\bar p$ be any perversity or superperversity. Let $\mc O$ be the orientation sheaf of the $n$-manifold $X-X^{n-1}$. Then $\mc P^*(\mc E\otimes \mc O)$ and $\mc{IS}^*(\mc E_0)$ are quasi-isomorphic, where $\mc{IS}^*(\mc E_0)$ is the  singular interesection chain sheaf with stratified coefficients.  
\end{theorem}

Of course one also modifies $\mc P^*$ in the obvious way so that the construction begins with 
 $\mc P^*_1=\mc E$ on $X-X^{n-1}$ and $$\mc P_{k+1}^*=\tau_{\leq \bar p(k)}Ri_{k*}\mc P_k^*$$ for $k\geq 1$. Also, in the proof, the axioms must be adjusted slightly in the obvious way to account for the codimension one stratum.

\section{Constructibility}\label{S: cc}

 Goresky and MacPherson initially built certain notions of sheaf constructibility into their axiomatic characterization of the Deligne sheaf on a pseudomanifold. Later, Borel showed in \cite[Section V.3]{Bo} that constructibility follows as a consequence of the other axioms. These arguments, however, use the local distinguished neighborhood structure of pseudomanifolds quite strongly, and thus it does not follow immediately from them that $\mc IS^*$ is constructible on an MHSS just because we have demonstrated that this sheaf complex satisfies the other axioms. So, in this section, we establish the desired constructibility properties of $\mc{IS^*}$. Since we showed in the last section that $\mc{IS^*}$ satisfies the axioms and since we observed that  any sheaf complex on an MHSS that satisfies the axioms is quasi-isomorphic to the Deligne sheaf, this implies the constructibility of any sheaf complex that satisfies the axioms on an MHSS. We note that, since pseudomanifolds \emph{are} MHSSs with sufficiently many approximate tubular neighborhoods (the distinguished neighborhoods), this also provides an alternative proof of constructibility of the Deligne sheaf on pseudomanifolds.

We first review the necessary concepts following the exposition in Borel \cite[Section V.3]{Bo}. All rings $R$ are noetherian commutative of finite cohomological dimension and possess a unity, and $X$ is locally compact of finite cohomological dimension over $R$. In particular, $X$ may be an $n$-dimensional MHSS - see Section \ref{S: MHSS}. 

\begin{definition}
\begin{itemize} 
\item  A direct system of $R$-modules $A_i$ is \emph{essentially constant} if for each $i$ in the index set $I$ there is an $i'\in I$, $i'\geq i$, such that $\ker(A_i\to A_{i'})=\ker(A_i\to \varinjlim A_j)$ and if there is an $i_0\in I$ such that $A_{i_0}\to \varinjlim A_j$ is surjective. 

\item  An  inverse system of $R$-modules $A_i$ is \emph{essentially constant} if for each $i$ in the index set $I$ there is an $i'\in I$, $i'\geq i$, such that $\im(A_{i'}\to A_{i})=\im(\varprojlim A_j\to A_i)$ and if there is an $i_0\in I$ such that $\varprojlim A_j\to A_{i_0}$ is injective. 
\end{itemize}
\end{definition}

If an inverse or direct system has a cofinal system, then it is essentially constant if and only if the cofinal system is. 

For the next definitions, we consider a bounded complex of sheaves $\mc S^*$ and a space $X$ filtered by closed subspaces. We let $\mc X$ denote the space together with its filtration information.

\begin{definition}
\begin{itemize}
\item $\mc S^*$ is \emph{cohomologically locally constant} (clc) if the derived sheaf $\mc H^*(\mc S^*)$ is a locally constant sheaf complex.

\item $\mc S^*$ is \emph{$\mc X$-cohomologically locally constant} ($\mc X$-clc) if $\mc H^*(\mc S^*)$ is locally constant on each stratum $X^i-X^{i-1}$.

\item $\mc S^*$ is \emph{$\mc X$-cohomologically constructible} ($\mc X$-cc) if it is $\mc X$-clc and, for each $x\in X$, the stalk $\mc H^*(\mc S^*)_x$ is finitely generated.

\item $\mc S^*$ is \emph{cohomologically constructible} (cc) if:
\begin{itemize}
\item For each $x\in X$ and $m\in \Z$, the inverse system $\H^m_c(U_x;\mc S^*)$ over all open neighborhoods of $x$ is essentially constant and its limit is finitely generated.

\item For each $x\in X$ and $m\in \Z$, the direct system $\H^m(U_x;\mc S^*)$ over all open neighborhoods of $x$ is essentially constant and its limit is finitely generated.

\item For each $x\in X$ and $m\in \Z$, $H^m(f^!_x\mc S^*)=\varprojlim \H^m_c(U_x;\mc S^*)$, where $U_x$ runs over open neighborhoods of $x$ and $f_x:x\into X$ is the inclusion. 

\item (Wilder's Property $(P,Q)$) If $P\subset Q$ are open in $X$, $\bar P\subset Q$, and $\bar P$ is compact, then the image of $\H^j_c(P;\mc S^*)$ in $\H^j(Q;\mc S^*)$ is finitely-generated for each $j$. 
\end{itemize}

\end{itemize}
\end{definition}

As observed in \cite[Section V.3.4]{Bo}, the four conditions for a sheaf complex to be cc are not independent; in fact the first condition implies the last two, and there are other interrelations. 

We will show that $\mc{IS}^*$ is $\mc X$-clc, $\mc X$-cc, and cc.

\begin{proposition}\label{P: clc}
Let $\mc{IS}^*$ be the intersection chain sheaf on the MHSS $X$ with sufficiently many local approximate tubular neighborhoods. Then $\mc{IS}^*$ is $\mc X$-clc. 
\end{proposition}
\begin{proof}
As seen in the proof of Proposition \ref{P: attaching range}, for any $x\in X$, $H^i(\mc{IS}^*_x)\cong IH^{\infty}_{n-i}(U)$ for any any local approximate tubular neighborhood $U$ of $x$. But if $x\in X_{n-k}$ and $y$ is another point in $U\cap X_{n-k}$, then $U$ is also a local approximate tubular neighborhood of $y$, and $H^i(\mc{IS}^*_y)\cong IH^{\infty}_{n-i}(U)\cong H^i(\mc{IS}^*_x)$. Since we also saw in the proof of Proposition \ref{P: attaching range} that the direct systems  $IH^{\infty}_{n-i}(V)$ are constant over cofinal sets of neighborhoods of $x$ and $y$, it follows that $H^i(\mc{IS}^*_x)$ is locally-constant over $X_{n-k}$. 
\end{proof}

\begin{theorem}
Let $\mc{IS}^*$ be the intersection chain sheaf on the MHSS $X$ with sufficiently many local approximate tubular neighborhoods. Then $\mc{IS}^*$ is $\mc X$-cc and cc. 
\end{theorem}
\begin{proof}
We have already seen that $\mc{IS}^*$ is $\mc X$-clc by Proposition \ref{P: clc}. Furthermore, we noted in the proof of that proposition that we have already seen in the proof of Proposition \ref{P: attaching range} that for any point $x\in X$, the direct system   $\H^*(U;\mc{IS}^*)\cong IH^{\infty}_{n-*}(U)$ is constant over a system of local approximate tubular neighborhoods  of $x$ that is cofinal in the direct system of all open neighborhoods of $x$. Thus 
$\H^*(V;\mc{IS}^*)$ is essentially constant over all open neighborhoods  $V$ of $X$ \cite[Remark V.3.2.b]{Bo}. 

Next we consider the inverse system $\H^*_c(V;\mc{IS}^*)$ over neighborhoods of $x\in X_{n-k}\subset X$. Recall once again that by \cite[Proposition 3.7]{GBF10}, the restriction of $\mc{IS^*}$ to any open set $V$ is quasi-isomorphic to the intersection chain sheaf of $V$. Since $X$ is locally-compact, the family of compact supports is paracompactifying, and by the same arguments as in the proof of Corollary 3.12 of \cite{GBF10}, $\H^*_c(V;\mc{IS}^*)\cong IH^c_{n-*}(V)$ (the cited Corollary assumes that $X$ is a topological stratified pseudomanifold, but this strict assumption is not necessary for the proof - the same arguments apply to any MHSS).  But now by \cite[Corollary 9.2 and Proposition 6.3]{GBF13},  if $U$ is a local approximate tubular neighborhood of $x$, then $IH^c_{n-*}(U)\cong IH^c_{n-*}(M_U)$, where $M_U$ is the mapping cylinder of the holink evaluation $\hl_s((X-X^{n-k})\cup (U\cap X_{n-k}), U\cap X_{n-k})\to U\cap X_{n-k}$. The mapping cylinder collapse $M_U\to U\cap X_{n-k}$ is a stratified fibration by \cite[Corollary 6.2]{Hug} and \cite[Proposition 3.3]{GBF5}. Furthermore, since $U\cap X_{n-k}$ is homeomorphic to $\R^{n-k}$, by \cite[Corollary 3.14]{GBF3}, $M_U$ is stratum- and fiber-preserving homotopy equivalent to $\R^{n-k}$ times the fiber over $x$, which is the cone on $\hl(X,x)$. The same is true then for any smaller local approximate tubular neighborhood $U'\subset U$ of $x$. Piecing together the appropriate stratified homotopy equivalences, one can see that the inclusion $IH^c_*(U')\to IH^c_*(U)$ is an isomorphism. (Alternatively, one could also use the long exact sequence of the pair and show that $IH^c_*(U,U')=0$ by using the MSAF property of approximate tubular neighborhoods to push any chain representing a relative cycle in $IC^c_*(U,U')$ into $IC^c_*(U')$.) It follows that the inverse system $\H^*_c(V;\mc{IS}^*)$ is also essentially constant. 

It remains to show for each $x\in X$ and a local approximate tubular neighborhood $U$ of $x$ that  $\H^*_c(U;\mc{IS}^*)\cong \varprojlim \H^*_c(V;\mc{IS}^*)$ and $\H^*(U;\mc{IS}^*)\cong \varinjlim\H^*(U;\mc{IS}^*)\cong H^*(\mc{IS}^*_x)$ are finitely generated. It will then follow from the definitions and \cite[Remarks V.3.4]{Bo} that $\mc{IS^*}$ is $\mc X$-cc and cc. 

We will proceed by induction over the strata of $X$. On the stratum $X-X^{n-2}$, $\mc{IS}^*$ is quasi-isomorphic to the sheaf of coefficents on a manifold and is both $\mc X$-cc and cc. Suppose now that $\mc{IS}^*|_{X-X^{n-k}}$ is $\mc X$-cc and cc, and let us add in the stratum $X_{n-k}$ and consider $\mc{IS}^*|_{X-X^{n-k-1}}$. Obviously, the local conditions that made $\mc{IS}^*|_{X-X^{n-k}}$ both $\mc X$-cc and cc continue to hold at points in $\mc{IS}^*|_{X-X^{n-k}}$, so we need only look at points in $X_{n-k}$ and show that the modules described in the last paragraph are finitely generated.

So let $x\in X_{n-k}$, $U$ a local approximate tubular neighborhood of $x$.  Once again, $\H^*_c(U;\mc{IS}^*)\cong IH^c_{n-*}(M)$, where $M$ is the mapping cylinder of the appropriate holink, and moreover $M$ is stratum-preserving homotopy equivalent to $\R^{n-k}\times c\ms L$, where $\ms L=\hl_s(X,x)$. By the cone formula then, $\H^*_c(U;\mc{IS}^*)$ is $0$ in some dimensions and isomorphic to $IH^c_{n-*}(\ms L)$ in others. Similarly, as calculated in the proof of Corollary \ref{C: ngbd reduce},  $\H^*(U;\mc{IS}^*)\cong IH^{c}_{n-*}(U, U-x)$, and by stratum-preserving homotopy equivalence, this is isomorphic to $IH^c_{n-*}(\R^{n-k}\times c\ms L, (\R^{n-k}\times c\ms L )-x)$. From the calculations of the proof of \cite[Proposition 2.20]{GBF10}, this is isomorphic to 
$IH^c_{k-*}(c\ms L, c\ms L-x)$. This too works out to be the compact intersection homology of $\ms L$ in some dimensions and $0$ in others. Thus it suffices to prove that $IH^c_*(\ms L)$, is finitely generated.

But it also follows from the various stratum-preserving homotopy equivalences we have employed that $\ms L$ must be stratum-preserving homotopy equialent to $U-U\cap X_{n-k}$. And we know $IH^c_*(U-U\cap X_{n-k})\cong \H_c^{n-*}(U-U\cap X_{n-k};\mc{IS}^*)$. Let $p: U-U\cap X_{n-k}\to \R^{n-k}\times \R\cong \R^{n-k+1}$ be the MSAF of the definition of an approximate tubular neighborhood. Let $D$ be the open unit disk in $\R^{n-k+1}$, and let $W$ be $p^{-1}(D)$. We will show that the inclusion $W\into U-U\cap X_{n-k}$, which induces $IH^c_*(W)\to IH^c_*(U-U\cap X_{n-k})$ and, equivalently, $\H_c^{n-*}(W;\mc{IS}^*)\to \H_c^{n-*}(U-U\cap X_{n-k};\mc{IS}^*)$ is an isomorphism. This will suffice, since by the induction hypothesis that $\mc{IS}^*|_{X-X^{n-k}}$ is cc, it must satisfy Wilder's $(P,Q)$ property. Here we take $Q=U-U\cap X_{n-k}$ and $P=W$, and we note that $\bar P\subset Q$ and $\bar P$ is compact (since $p$ is proper). The Wilder property then allows us to conclude that the image of $\H_c^{n-*}(W)$ in  $\H_c^{n-*}(U-U\cap X_{n-k})$, which is equal to $\H_c^{n-*}(U-U\cap X_{n-k})$, is finitely generated. 

So now to complete the proof, consider the exact sequence of the pair $(U-U\cap X_{n-k},W)$.  We show that $IH^c_*(U-U\cap X_{n-k},W)=0$. Let $\xi$ be a relative cylce in $IC^c_*(U-U\cap X_{n-k},W)=0$. Let $r:\R^{n-k+1}\times I\to \R^{n-k+1}$ be a radial deformation retraction from the identity map to the collapse map to the origin. Consider $F=r(p\times \text{id}_I):|\xi|\times I\to \R^{n-k+1}$. Since $p$ is an MSAF, there is a stratified approximate lift $\td F: |\xi|\times I\times [0,1)\to U-U\cap X_{n-k}$, and the associated map $\bar F: |\xi|\times I\times I\to \R^{n-k}$. Note that $\bar F(|\xi|\times 1\times 1)=0\subset \R^{n-k+1}$, so $|\xi|\times  1\times 1\subset \bar F^{-1}(D)$, which is an open set. Since $|\xi|$ is compact, it follows from elementary topology, that there is an open neighborhood $A$ of $1\times 1$ in $I\times I$ such that $|\xi|\times A\subset \bar F^{-1}(D)$. Similarly, since $\bar F(|\bd \xi|\times I\times 1)\subset D$ and $|\bd \xi|\times I$ is compact, there is a neighborhoof $B$ of $1$ in $I$ such that $|\bd \xi|\times I\times B\subset \bar F^{-1}(D)$. So now we choose a path $\gamma$ in $I\times [0,1)$ such that 
\begin{enumerate}
\item $\gamma(0)\in 0\times [0,1)$, 
\item $\gamma(1)\in 1\times [0,1)$, 
\item $\gamma \subset I\times B$, and 
\item $\gamma(1)\in A$.
\end{enumerate}
Then the homotopy $\td F\circ (\text{id}_{|\xi|}\times \gamma): |\xi|\times I\to U-U_{n-k}$ retracts $|\xi|$ into $W$ in a stratum-preserving manner and keeps $\bd \xi$ in $W$. Thus this homotopy can be used to construct a relative null-homology. 
\end{proof}

\section{Poincar\'e Duality}\label{S: pd}

The initial impetus for the study of intersection homology was the goal of extending Poincar\'e duality to manifold stratified spaces. This was first achieved with field coefficients for compact PL pseudomanifolds by Goresky and MacPherson in \cite{GM1}, where it was shown that if $X$ is an $n$-dimensional compact oriented PL stratified pseudomanifold and if $\bar p$ and $\bar q$ are dual perversities ($\bar p(k)+\bar q(k)=k-2$), then there is a nonsingular pairing $I^{\bar p}H_i(X;\Q)\otimes  I^{\bar q}H_{n-i}(X;\Q)\to \Q$. If $X$ has only even-codimension singularities (or more generally if $X$ is a \emph{Witt space} - see \cite{Si83}), then the \emph{upper and lower middle perversities}, $(0,1,1,2,2,3,\ldots)$ and $(0,0,1,1,2,2,\ldots)$, yield isomorphic intersection homology groups, and there is a pairing $I^{\bar m}H_i(X;\Q)\otimes  I^{\bar m}H_{n-i}(X;\Q)\to \Q$. If, in addition, $n=2k$, one obtains an $\epsilon$-symmetric self-pairing
$I^{\bar m}H_k(X;\Q)\otimes  I^{\bar m}H_{k}(X;\Q)\to \Q$, which leads to signature invariants, $L$-classes, etc. 
Using sheaf-theoretic machinery, this version of Poincar\'e duality and its consequences were extended to topological pseudomanifolds  and more general coefficient systems over fields in \cite{GM2} (see also \cite[Section V.9]{Bo}). 

Goresky and Siegel then showed in \cite{GS83} that Poincar\'e duality on pseudomanifolds holds over the integers, provided  certain torsion subgroups of the intersection homology groups of all links vanishes. In particular, they defined a pseudomanifold to be \emph{locally $\bar p$-torsion free} if, for all $k$ and for each $x\in X_{n-k}$ with corresponding link $L_x$, $I^{\bar p}H_{k-2-\bar p(k)}(L_x)$ is torsion free. With this assumption, one obtains a nonsingular pairing
$$I^{\bar p}H_i(X)/\text{torsion}\otimes  I^{\bar q}H_{n-i}(X)/\text{torsion}\to \Z,$$ as well as a nonsingular torsion pairing
$$T^{\bar p}H_i(X)\otimes  T^{\bar q}H_{n-i-1}(X)\to \Q/\Z,$$
where $TH$ denotes the torsion subgroup of $IH$.  

We now show that this version of Poincar\'e duality further extends to include MHSSs with sufficiently many local approximate tubular neighborhoods. It will follow from a theorem of Quinn concerning the topological invariance of $IH^c_*(X)$ for MHSSs that $IH^{\infty}_*(X)$ is also a topological invariant, assuming sufficiently many approximate tubular neighborhoods.

First, we need an analogue of the Goresky-Siegel link condition:

\begin{definition}
Let $R$ be a PID. We say that the MHSS $X$ is \emph{homotopy locally ($\bar p$,$R$)-torsion free} if  for all $k$ and each $x\in X_{n-k}$,  $I^{\bar p}H^c_{k-2-\bar p(k)}(\mc L_x)$ is $R$-torsion free, where $\mc L_x$ is the homotopy link of $x$ in $X$. Utilizing the computations as in the proof of Proposition \ref{P: 0 range}, this condition is equivalent to $H_c^{n-k+2+\bar p(k)}(U;\mc{I}^{\bar p}\mc S^*)$ being torsion free for any local approximate tubular neighborhood $U$ of $x$.

 This definition is a direct analogue of the definition of locally $\bar p$-torsion free in Goresky-Siegel \cite{GS83}.
Note that any $X$ is automatically homotopy locally ($\bar p$,$R$)-torsion free if $R$ is a field.
\end{definition}

\begin{comment}
\begin{definition}
Let $R$ be a PID. We say that the MHSS $X$ is \emph{homotopy locally ($\bar p$,$R$)-torsion free} if  for all $k$ and each $x\in X_{n-k}$,  $H_c^{n-k+2+\bar p(k)}(\mc{I}^{\bar p}\mc S^*_x)$ is $R$-torsion free.   
Utilizing the computations as in the proof of Proposition \ref{P: 0 range}, $H_c^{n-k+2+\bar p(k)}(\mc{I}^{\bar p}\mc S^*_x)\cong I^{\bar p}H^c_{k-2-\bar p(k)}(\mc L_x)$, where $\mc L_x$ is the homotopy link of $x$ in $X$. Thus this definition is a direct analogue of the definition of locally $\bar p$-torsion free in Goresky-Siegel \cite{GS83}.
Note that any $X$ is automatically homotopy locally ($\bar p$,$R$)-torsion free if $R$ is a field.
\end{definition}

\end{comment}

This leads to our main theorem:
\begin{comment}
\begin{theorem}\label{T: PD}
Let $X$ be an n-dimensional MHSS with no codimension one stratum  and with sufficiently many local approximate tubular neighborhoods. Let $\mc O$ be the orientation sheaf of the $n$-manifold $X-X^{n-2}$, and let $\mc E$ be a local coefficient system of finitely-generated $\F$-modules for a field $\F$. Let $\bar p$ and $\bar q$ be dual perversities ($\bar p(k)+\bar q(k)=k-2$). Then $$I^{\bar q}H^{\infty}_i(X;\Hom(\mc E, \F_{X-X^{n-2}})\otimes \mc O)\cong \text{\emph{Hom}}(I^{\bar p}H^c_{n-i}(X;\mc E);\F).$$
\end{theorem}

\end{comment}

\begin{theorem}\label{T: PD}
Let $X$ be a homotopy locally ($\bar p$,$R$)-torsion free n-dimensional MHSS with no codimension one stratum  and with sufficiently many local approximate tubular neighborhoods. Let $\mc O$ be the orientation sheaf of the $n$-manifold $X-X^{n-2}$, and let $\mc E$ be a local coefficient system on $X-X^{n-2}$ of finitely-generated free modules over the  principal ideal domain $R$. Let $\bar p$ and $\bar q$ be dual perversities ($\bar p(k)+\bar q(k)=k-2$).
Let $TH_*$ and $FH_*$ denote, respectively, the $R$-torsion subgroup and $R$-torsion free quotient group of $IH_*$, and let $Q(R)$ denote the field of fractions of $R$.

Suppose that $\Hom(T^{\bar p}H^c_{i-1}(X;\mc E),Q(R)/R)$ is a torsion $R$-module (in particular, if $T^{\bar p}H^c_{i-1}(X;\mc E)$ is finitely generated).
Then $$\Hom(F^{\bar p}H^c_i(X;\mc E),R)\cong F^{\bar q}H^{\infty}_{n-i}(X;\Hom(\mc E, R_{X-X^{n-2}})\otimes \mc O)$$ and $$\Hom(T^{\bar p}H^c_{i-1}(X;\mc E),Q(R)/R)\cong T^{\bar q}H^{\infty}_{n-i}(X;\Hom(\mc E, R_{X-X^{n-2}})\otimes \mc O).$$  
\end{theorem}

We record separately the case for field coefficients, for which all of the torsion conditions are satisfied automatically.

\begin{corollary}\label{C: PD}
Let $X$ be an n-dimensional MHSS with no codimension one stratum  and with sufficiently many local approximate tubular neighborhoods. Let $\mc O$ be the orientation sheaf of the $n$-manifold $X-X^{n-2}$, and let $\mc E$ be a local coefficient system on $X-X^{n-2}$ of finitely-generated $\F$-modules for a field $\F$. Let $\bar p$ and $\bar q$ be dual perversities ($\bar p(k)+\bar q(k)=k-2$). Then $$\text{\emph{Hom}}(I^{\bar p}H^c_{n-i}(X;\mc E);\F)\cong I^{\bar q}H^{\infty}_i(X;\Hom(\mc E, \F_{X-X^{n-2}})\otimes \mc O).$$
\end{corollary}

When $X$ is compact and orientable, we obtain as a special case the simpler, but more familiar, statement $$\Hom(I^{\bar p}H_i(X;\Q),\Q)\cong I^{\bar q}H_{n-i}(X;\Q).$$ If, in addition, $X$ is homotopy locally ($\bar p$,$\Z$)-torsion free, we have
$$\Hom(F^{\bar p}H_i(X),\Z)\cong F^{\bar q}H_{n-i}(X) \qquad \text{and} \qquad \Hom(T^{\bar p}H_i(X),\Q/\Z)\cong T^{\bar q}H_{n-i}(X).$$

\begin{proof}[Proof of Theorem \ref{T: PD}]
Having established in previous sections that the singular chain intersection homology on $X$ corresponds to the Deligne sheaf hypercohomology and having used this correspondence to establish the constructibility of $\mc IS^*\sim_{q.i.} \mc P$, the main idea of the proof of Poincar\'e duality is the same as that in prior treatments for pseudomanifolds: We consider the Verdier dual of an intersection chain sheaf and show that it satisfies the axioms that make it the intersection chain sheaf with the dual perversity.\footnote{In  \cite{GM2} and \cite{Bo}, it is shown for topological pseudomanifolds that the dual of the Deligne sheaf satisfies the axioms AX2. We instead  continue to utilize the axioms AX1, as already presented above.} Many of the details, though, rely on the properties we have divined for the sheaf of singuler intersection chains. 

Let $\mc D^*_X$ be the Verdier dualizing functor on $X$, which takes a bounded sheaf complex $\mc A^*$ in $D^b(X)$, the bounded derived category of sheaves on $X$,  to $\text{\emph{RHom}}^*(\mc A^*,\D_X^*)$, where $\D_X^*=f^!R_{pt}$, $f$ the map from $X$ to a point. 
Thorough expositions on $\D_X^*$ and the functor $f^!$ can be found in both \cite{Bo} and \cite{BaIH}. Below, we will show the following:

\begin{lemma}\label{qi} Under the assumptions of the theorem, $\mc D^*_X(\mc I^{\bar p}\mc S^*(\mc E))[-n]$ is quasi-isomorphic to $\mc I^{\bar q}\mc S^*(\mc D^*_{X-X^{n-2}}(\mc E)[-n])$, where $[-n]$ is the shift such that for a complex $\mc A^*$, $(\mc A^*[-n])^i=\mc A^{i-n}$.
\end{lemma}

From this lemma, the duality theorem is proven just as in Goresky-Siegel \cite{GS83} as follows:

 Given $\mc D^*(\mc S^*)$, for any sheaf $\mc S^*$ over a principal ideal domain, there is a short exact sequence (see \cite[V.7]{Bo})
\begin{diagram}
0&\rTo & \Ext(\H^{i+1}_c(X;\mc S^*),R)  &\rTo & \H^{-*}(X;\mc D^*(\mc S^*)) &\rTo &\Hom(\H^i_c(X;\mc S^*),R) &\rTo & 0
\end{diagram}

Applying the lemma with $\mc S^*=\mc I^{\bar p}\mc S^*(\mc E)$, we have 
\begin{align*}
\H^{-*}(X;\mc D^*\mc I^{\bar p}\mc S^*(\mc E))& \cong \H^{n-*}(X;\mc D^*\mc I^{\bar p}\mc S^*(\mc E)[-n])\\
&\cong \H^{n-*}(X;\mc I^{\bar q}\mc S^*(\mc D^*_{X-X^{n-2}}(\mc E)[-n]))\\
&\cong I^{\bar q}H^{\infty}_*(X;\mc D^*_{X-X^{n-2}}(\mc E)[-n]).
\end{align*} 

Recall that, for a principal ideal domain, $Ext(A,R)\cong \Hom(T(A),Q(R)/R)$, where $T(A)$ is the $R$-torsion subgroup of $A$ and $Q(R)$ is the field of fractions of $R$. 
 So the preceding exact sequence becomes

\begin{diagram}
0&\rTo & \Hom(T^{\bar p}H^c_{i-1}(X;\mc E),Q(R)/R)  &\rTo & I^{\bar q}H^{\infty}_*(X;\mc D^*(\mc E)[-n]) &\rTo &\Hom(I^{\bar p}H^c_i(X;\mc E),R) &\rTo & 0.
\end{diagram}

Since $\Hom(T^{\bar p}H^c_{i-1}(X;\mc E),Q(R)/R)$ is $R$-torsion by hypothesis and $\Hom(I^{\bar p}H^c_i(X;\mc E),R)$ must be torsion free, this exact sequence is naturally isomorphic to 

\begin{diagram}
0&\rTo &  T^{\bar q}H^{\infty}_*(X;\mc D^*(\mc E)[-n]) &\rTo & I^{\bar q}H^{\infty}_*(X;\mc D^*(\mc E)[-n])  &\rTo &F^{\bar q}H^{\infty}_*(X;\mc D^*(\mc E)[-n]) &\rTo & 0,
\end{diagram}
where the first nontrivial map is simply the inclusion of the torsion subgroup.

Thus we obtain isomorphisms $T^{\bar q}H^{\infty}_*(X;\mc D^*(\mc E^*)[-n])\cong \Hom(T^{\bar p}H^c_{i-1}(X;\mc E),Q(R)/R)$ and $F^{\bar q}H^{\infty}_*(X;\mc D^*(\mc E^*)[-n])\cong \Hom(I^{\bar p}H^c_i(X;\mc E),R)$.

Finally, we note that $\mc D^*(\mc E)[-n]\cong \Hom(\mc E, R_{X-X^{n-2}})\otimes \mc O$ by \cite[V.7.10.4]{Bo}.
\end{proof}

We now prove the above-stated Lemma \ref{qi}, showing that the Verdier dual of a perversity $\bar p$ intersection chain sheaf is a perversity $\bar q$ intersection chain sheaf.

\begin{proof}[Proof of Lemma  \ref{qi}]
As noted in the proof of Theorem \ref{T: IS is deligne}, by \cite[Theorem V.2.5]{Bo} it suffices to show that 
$\mc D^*(\mc I^{\bar p}\mc S^*(\mc E))[-n]$ satisfies the axioms $AX1_{\bar q, X}(\mc D^*_{X-X^{n-2}}(\mc E)[-n])$.

Let $X-X^{n-2}=U_2$. By \cite[V.7.10(4)]{Bo}, $\mc D^*_{U_2}(\mc E)[-n]\sim_{q.i} \text{\emph{Hom}}(\mc E, R_{U_2})\otimes \mc O$, so $\mc D^*_{U_2}(\mc E)[-n]$ is indeed a local system of coefficients. 
Then, $(\mc D^*(\mc I^{\bar p}\mc S^*(\mc E))[-n])|_{U_2}=\mc D_{U_2}^*(\mc I^{\bar p}\mc S^*(\mc E)|_{U_2})[-n]\cong \mc D_{U_2}^*(\mc E)[-n]$ by \cite[VI.3.11.2]{Bo}, \cite[V.10.11]{Bo}, and Axiom $AX1_{\bar p}(\mc E)a$ for $\mc I^{\bar p}\mc S^*(\mc E)$. This establishes the last part of axiom 1.

Next, let $x\in X^{n-k}-X^{n-k-1}$ and consider $H^*(\mc D^*(\mc I^{\bar p}\mc S^*(\mc E))_x[-n])\cong \varinjlim_{x\in U} \H^{*-n}(U;\mc D^*(\mc I^{\bar p}\mc S^*(\mc E)))$. For any sheaf complex $\mc A^*$ over $R$ in $D^b(X)$ and any open set $U\in X$, we have an exact sequence
 $$0\to \Ext(\H_c^{i+1}(U;\mc A^*),R)\to \H^{-i}(U;\mc D_X \mc A^*)\to \Hom(\H_c^i(U;\mc A^*),R)\to 0$$ (see \cite[Section 3.4]{BaIH}). Thus there is an exact sequence
 $$0\to\Ext(IH^c_{*-1}(U;\mc E),R) \to \H^{*-n}(U;\mc D^*(\mc I^{\bar p}\mc S^*(\mc E)))\to \Hom(IH^c_{*}(U;\mc E),R)\to 0.$$ As shown in the proof of Proposition \ref{P: 0 range}, if $U$ is a local approxiate tubular neighborhood, then $U$ is stratum-preserving homotopy equivalent to a cone $c\ms L$, where we can treat $\ms L$ as a $k-1$ dimensional filtered space. Thus by the standard cone formula for singular intersection homology \cite{Ki}, $IH^c_{*}(U;\mc E)=0$  for $*> k-2-\bar p(k)=k-2-(k-2-\bar q(k))=\bar q(k)$. Furthermore, since $X$ is homotopy locally ($\bar p$,$R$)-torsion free, $\Ext(IH^c_{k-2-\bar p(k)}(U;\mc E),R)$ is also $0$, so  $\H^{*-n}(U;\mc D^*(\mc I^{\bar p}\mc S^*(\mc E)))=0$ for $*>\bar q(k)$. It is also clear that these groups must be $0$ for $*<0$.

Finally, to verify the attaching axiom, we observe (by the discussion following the statement of the axioms, above) that for any sheaf $\mc A^*$, the attaching map induces an isomorphism of cohomology stalks at $x\in X_{n-k}=X^{n-k}-X^{n-k-1}$ in dimension $j$ if and only if restriction induces an isomorphism $\varinjlim_{x\in U} \H^j(U;\mc A^*)\to \varinjlim_{x\in U}\H^j(U-U\cap X_{n-k};\mc A^*)$, where $U$ runs over all open neighborhoods of $x$. Of course, we can limit ourselves to a cofinal system of local approximate tubular neighborhoods, and it suffices to find then isomorphisms $\H^j(U;\mc A^*)\to \H^j(U-U\cap X_{n-k};\mc A^*)$ that are functorial in that they commute with further restrictions. In the case at hand, we study $\H^j(U;\mc D^*(\mc I^{\bar p}\mc S^*(\mc E))[-n])\to \H^j(U-U\cap X_{n-k};\mc D^*(\mc I^{\bar p}\mc S^*(\mc E))[-n])$, which induces a map of short exact sequences

\begin{diagram}
0&&0\\
\dTo&&\dTo\\
 \Ext(\H^{n-j+1}_c(U;\mc I^{\bar p}\mc S(\mc E)),R)&\rTo& \Ext(\H^{n-j+1}_c(U-U\cap X_{n-k};\mc I^{\bar p}\mc S^*(\mc E)),R) \\
 \dTo&&\dTo\\
 \H^j(U;\mc D^*(\mc I^{\bar p}\mc S^*(\mc E))[-n]) &\rTo&\H^j(U-U\cap X_{n-k};\mc D^*(\mc I^{\bar p}\mc S^*(\mc E))[-n])\\
 \dTo&&\dTo\\
\Hom(\H^{n-j}_c(U;\mc I^{\bar p}\mc S(\mc E)),R) &\rTo&\Hom(\H^{n-j}_c(U-U\cap X_{n-k};\mc I^{\bar p}\mc S^*(\mc E)),R) \\
 \dTo&&\dTo\\
0&&0&,\\
\end{diagram}

\begin{comment}
\begin{diagram}
0&\rTo & \Ext(\H^{n-j+1}_c(U;\mc I^{\bar p}\mc S(\mc E)),R)  &\rTo & \H^j(U;\mc D^*(\mc I^{\bar p}\mc S^*(\mc E))[-n])  &\rTo & \Hom(\H^{n-j}_c(U;\mc I^{\bar p}\mc S(\mc E)),R)  &\rTo & 0\\
 &     &\dTo & &\dTo & &\dTo &\\
0&\rTo & \Ext(\H^{n-j+1}_c(U-U\cap X_{n-k};\mc I^{\bar p}\mc S^*(\mc E)),R)   &\rTo & \H^j(U-U\cap X_{n-k};\mc D^*(\mc I^{\bar p}\mc S^*(\mc E))[-n])  &\rTo & \Hom(\H^{n-j}_c(U-U\cap X_{n-k};\mc I^{\bar p}\mc S^*(\mc E)),R)  &\rTo & 0,
\end{diagram}
\end{comment}
\begin{comment}
functorially isomorphic  to $\Hom(\H^{n-j}_c(U;\mc I^{\bar p}\mc S(\mc E)),\F)\to \Hom(\H^{n-j}_c(U-U\cap X_{n-k};\mc I^{\bar p}\mc S^*(\mc E)),\F)$ or, equivalently, $\Hom(I^{\bar p}H^c_{j}(U;\mc E),\F)\to \Hom(I^{\bar p}H^c_{j}(U-U\cap X_{n-k};\mc E),\F)$,
\end{comment}

\noindent where the maps of the outer terms are induced by the inclusion maps $I^{\bar p}H^c_{j}(U-U\cap X_{n-k};\mc E)\to I^{\bar p}H^c_{j}(U;\mc E)$; we present a proof of this below in  Appendix A. Once again, we know that $\H^{n-*}_c(U-U\cap X_{n-k};\mc I^{\bar p}\mc S(\mc E))\cong I^{\bar p}H^c_{*}(U-U\cap X_{n-k};\mc E)$ is isomorphic to $I^{\bar p}H^c_*(\ms L;\mc E)$ and $\H^{n-*}_c(U;\mc I^{\bar p}\mc S(\mc E))\cong I^{\bar p}H^c_{*}(U;\mc E)\cong I^{\bar p}H^c_*(c\ms L;\mc E)$, where $\ms L\sim \hl(X,x)$. By the cone formula, the inclusion $I^{\bar p}H_*^c(\ms L;\mc E)\to I^{\bar p}H_*^c(c\ms L;\mc E)$  is an isomorphism for $*<k-1-\bar p(k)$. Thus, by the five lemma, $\H^j(U;\mc I^{\bar p}\mc S(\mc E)[-n])\to \H^j(U-U\cap X_{n-k};\mc I^{\bar p}\mc S^*(\mc E))[-n]$ is an isomorphism for $j\leq \bar q(k)$. Since this computation is functorial with respect to restrictions, we obtain the desired isomorphism in the limits.

Thus $\mc D^*(\mc I^{\bar p}\mc S^*(\mc E))[-n]$ satisfies the axioms $AX1_{\bar q, X}(\mc D^*_{X-X^{n-2}}(\mc E)[-n])$, which completes the proof of the lemma. 
\end{proof}

\begin{comment}
using the well-known arguments that follow from the lemma:
\begin{align*}
I^{\bar q}H^{\infty}_i(X;\text{\emph{Hom}}(\mc E, \F_{X-X^{n-2}})\otimes \mc O)&\cong \H^{n-i}(X;\mc{I}^{\bar q}\mc S^*(\text{\emph{Hom}}(\mc E, \F_{X-X^{n-2}})\otimes \mc O))\\
&\cong \H^{n-i}(X;\mc{I}^{\bar q}\mc S^*(\mc D^*_{X-X^{n-2}}(\mc E)[-n]))\\
&\cong \H^{n-i}(X;\mc D^*(\mc I^{\bar p}\mc S(\mc E))[-n])\\
&\cong \Hom(\H^i_c(X;\mc I^{\bar p}\mc S(\mc E));\F)\\
&\cong \Hom(IH_{n-i}^c(X; \mc E);\F).
\end{align*}
In going from the first to the second line, we have used that $\text{\emph{Hom}}(\mc E, \F_{X-X^{n-2}})\otimes \mc O\cong \mc D^*_{X-X^{n-2}}(\mc E)[-n]$; see \cite[Section V.7.10]{Bo}.
\end{comment}

\begin{corollary}\label{C: invariance}
Let $X$ be a  homotopy locally ($\bar p$,$R$)-torsion free n-dimensional MHSS with no codimension one stratum  and with sufficiently many local approximate tubular neighborhoods. Suppose $\mc E$ is a local coefficient system on $X-X^{n-2}$ of finitely-generated free $R$ modules over a PID $R$. Then $I^{\bar q}H^{\infty}_*(X;\mc E)$ is a topological invariant, i.e. it does not depend on the choice of stratification of $X$ as an MHSS.
\end{corollary}
\begin{proof}
Let $\mc E^*=\text{\emph{Hom}}(\mc F;R_{X-X^{n-2}})\otimes \mc O$. Then $\mc E=\text{\emph{Hom}}(\mc E^*\otimes \mc O;R)$. By the proof of the theorem, $I^{\bar q}H^{\infty}_*(X;\mc E)$ is part of a short exact sequence with  $\Hom(I^{\bar p}H^c_{n-*}(X;\mc E^*),R)$ and $\Ext(I^{\bar p}H^c_{n-*-1}(X;\mc E^*),R)$. But according to \cite[Section 2]{Q2}, $I^{\bar p}H^c_{*}(X;\mc E^*)$ is independent of the stratification of $X$, and thus the same must follow for $I^{\bar q}H^{\infty}_*(X;\mc E)$ by the five lemma.
\end{proof}

\section{Homotopy Witt spaces and Poincar\'e Duality Spaces}\label{S: apps}

Let $\bar m$ and $\bar n$ be the lower-middle and upper-middle perversities:
$\bar m(k)=\lfloor \frac{k-2}{2} \rfloor$ and $\bar n(k)=\lfloor \frac{k-1}{2} \rfloor$. Let $R$ be a fixed principal ideal domain.

Generalizing the definition of Siegel \cite{Si83}, we can define \emph{homotopy Witt spaces}:

\begin{definition}
We say that the compact homotopy locally ($\bar m$,$R$)-torsion free $n$-dimensional  MHSS with sufficiently many approximate tubular neighborhoods $X$ is a  \emph{homotopy $R$-Witt space} (or $HR$-Witt space) if for each $x$ in each odd-codimension stratum $X_{n-(2k+1)}$, we have $H^{\bar n(k)}(\mc{I}^{\bar n}\mc S^*_x;R)=H^{k}(\mc{I}^{\bar n}\mc S^*_x;R)=0$.  
\end{definition}
Utilizing the computations as in the proof of Proposition \ref{P: 0 range}, $H^{k}(\mc{I}^{\bar n}\mc S^*(R)_x)\cong I^{\bar n}H^c_{k}(\mc L;R)$, where $\mc L$ is the homotopy link of $x$ in $X$. \begin{comment}We observe that for intersection homology the universal coefficients theorem does hold between $\Z$ and $\Q$ by \cite[Theorem 8.1]{GS83} so that this is  $I^{\bar n}H^c_{k}(\mc L;\Q)$. Thus an equivalent definition for homotopy rational Witt spaces would be that $I^{\bar n}H^c_{k}(\mc L;\Q)=0$ for homotopy links of points in strata of dimension $2k+1$. This is more closely analagous to the original definition of Siegel.\end{comment} 
We note that there is a slight difference from the usual definition of Witt spaces in that our formula uses the upper-middle perversity. For Siegel this is not an issue because for $L$ a compact orientable pseudomanifold, $I^{\bar n}H^c_{k}(L;\Q)\cong I^{\bar m}H^c_{k}(L;\Q)$ by the intersection homology Poincar\'e duality of Goresky-MacPherson. But we cannot assume that we have such an isomorphism for the homotopy link $\mc L$. 

We note, incidentally, that our ``Witt condition'' implies that $X$ has no codimension one stratum since, if $x\in X_{n-1}$, then  $H^{0}(\mc{I}^{\bar n}\mc S^*(R)_x)\cong H^0((Ri_{1*}R_{X-X^{n-1}})_x)\cong \dlim_{x\in U}\H^0(U-U\cap X^{n-1};R)\cong \dlim_{x\in U}H^0(U-U\cap X^{n-1};R)$.  This will never be $0$.

If $X$ is an $HR$-Witt space, then, by Theorem \ref{T: IS is deligne}, each $\mc{I}^{\bar p}\mc S^*(R)$ is quasi-isomorphic to the Deligne sheaf with the appropriate perversity, and it then follow immediately, as for Siegel's Witt spaces, that $\mc{I}^{\bar m}\mc S^*(R)$ is quasi-isomorphic to $\mc{I}^{\bar n}\mc S^*(R)$: each inclusion $\tau_{\leq \bar m(k)}Ri_{k*}\mc P_k^*\into \tau_{\leq \bar n(k)}Ri_{k*}\mc P_k^*$ is a quasi-isomorphism. By \cite[Section 2]{Q2}, these groups are topological invariants, so in fact we will have $\mc{I}^{\bar m}\mc S^*(R)\sim_{q.i.}\mc{I}^{\bar n}\mc S^*(R)$ if the topological space $X$ can be given the structure of an $HR$-Witt space with respect to any stratification.

If $X$ is a compact orientable $HR$-Witt space, then we have from Theorem \ref{T: PD} that $F^{\bar m}H_i(X;R)\cong \text{\emph{Hom}}(F^{\bar m}H_{n-i}(X;R);R)$ and $T^{\bar m}H_i(X;R)\cong \text{\emph{Hom}}(T^{\bar m}H_{n-i-1}(X;R);Q(R)/R)$.
In particular, we have the following theorem

\begin{theorem} If $X$ is a compact orientable $HR$-Witt space of dimension $2n$, there is a nonsingular pairing $F^{\bar m}H_n(X;R)\otimes F^{\bar m}H_{n}(X;R)\to R$. If $n$ is even, then $X$ has a well-defined signature. 
\end{theorem}

\begin{comment}
We can then ask the following natural question: how much of the machinery of characteristic classes and self-dual sheaves on PL Witt spaces and other topological pseudomanifolds carries over to $H\Q$-Witt spaces and other HMSSs?
                                                                                                \end{comment}

\section{More general ground rings}\label{S: rings}

As shown by Goresky and MacPherson for pseudomanifolds \cite{GM2}, the duality quasi-isomorphism  $\mc D^*_X(\mc I^{\bar p}\mc S^*(\mc E))[-n]\sim_{q.i.} \mc I^{\bar q}\mc S^*(\mc D^*_{X-X^{n-2}}(\mc E)[-n])$ holds with field coefficients with no further assumptions on the properties of the space. Goresky and Siegel \cite{GS83} extended this duality to integer coefficients (though their argument would work for any principal ideal domain) at the expense of requiring a single torsion condition on the links of each stratum. Above, we have considered the analogous conditions and duality quasi-isomorphisms on MHSSs. In this section, we explore what conditions may be imposed on our space in order for $\mc D^*_X(\mc I^{\bar p}\mc S^*(\mc E))[-n]\sim_{q.i.} \mc I^{\bar q}\mc S^*(\mc D^*_{X-X^{n-2}}(\mc E)[-n])$ to hold for rings of greater cohomological dimension. 

Throughout this section, let $R$ be a fixed noetherian commutative ring of finite cohomological dimension.

We must examine the proof of Lemma \ref{qi}, in which we demonstrated that $\mc D^*_X(\mc I^{\bar p}\mc S^*(\mc E))[-n]$ satisfies the Goresky-MacPherson axioms $AX1_{\bar q, X}(\mc D^*_{X-X^{n-2}}(\mc E)[-n])$. There is no problem with showing for any $R$ that $\mc D^*_{X-X^{n-2}}(\mc I^{\bar p}\mc S^*(\mc E)|_{X-X^{n-2}})[-n]\sim_{q.i.}\mc D^*_{X-X^{n-2}}(\mc E)[-n]$ in the same manner as above, so we move on to the other axioms. 

In the proof of Lemma \ref{qi}, we used the ``universal coefficient'' short exact sequence to show that $H^*(\mc D^*(\mc I^{\bar p}\mc S^*(\mc E))_x[-n])=0$ for $*>\bar q(k)$. If $R$ is not a PID, we will not have this exact sequence in general, but we will have a spectral sequence instead. In general, for any sheaf complex $\mc A^*$, there is a spectral sequence abutting to $\H^*(U;\mc D^*\mc A^*)$ with $E_2$ terms $E_2^{r,s}\cong \Ext^r(\H^{-s}_c(U;\mc A^*),R)$ (see \cite[Section V.7.7]{Bo}). When $R$ is a PID, it is the collapsing of this spectral sequence at the $E_2$ terms that leads to the short exact sequence. 

Now, suppose $x\in X_{n-k}$ and let $U$ be a local approximate tubular neighborhood of $x$.
At the $\infty$ stage, the terms that will influence  $\H^i(U;\mc D^*(\mc I^{\bar p}\mc S^*(\mc E))[-n])\cong \H^{i-n}(U;\mc D^*(\mc I^{\bar p}\mc S^*(\mc E)))$ are  the terms $E_{\infty}^{r,i-n-r}\cong \Ext^r(\H^{r+n-i}_c(U;\mc I^{\bar p}\mc S^*(\mc E)),R)\cong \Ext^r(I^{\bar p}H_{i-r}^c(U;\mc E),R)$. So, a sufficient condition to guarantee that $H^*(\mc D^*(\mc I^{\bar p}\mc S^*(\mc E))_x[-n])=0$ for $*>\bar q(k)$ would be to ask that $\Ext^r(I^{\bar p}H_{i-r}^c(U;\mc E))$ vanishes for $i>\bar q(k)$ and for all $r$. Recall that  $\bar q(k)=k-2-\bar p(k)$ and that $U\sim_{s.p.h.e}c\mc L_x\times R^{n-k}$, where $\mc L_x$ has formal dimension $k-1$. Thus when $r=0$, this condition is satisfied automatically due to the usual intersection homology cone formula, according to which $I^{\bar p}H_{i}^c(cL;\mc E)=0$ for $i>k-2-\bar p(k)$, and for $r=1$, this is precisely our earlier homotopy locally ($\bar p$,$R$)-torsion free condition, which generalizes the Goresky-Siegel condition.

The situation for the attaching axiom is a little more complicated. We continue to let $U$ be a local approximate tubular neighborhood of $x\in X_{n-k}$. We would like for the restrictions
$\H^i(U;\mc D^*(\mc I^{\bar p}\mc S^*(\mc E))[-n])\to \H^i(U-U\cap X_{n-k};\mc D^*(\mc I^{\bar p}\mc S^*(\mc E))[-n])$ to be an isomorphism for $i\leq \bar q(k)$. To study this map, we turn to the corresponding map of spectral sequences, which we will denote $E(U)\to E(U-U\cap X^{n-k})$. The $E_2$ term maps will be $\Ext^r(I^{\bar p}H_{i-r}^c(U;\mc E),R)\to \Ext^r(I^{\bar p}H_{i-r}^c(U-U\cap X^{n-k};\mc E),R)$. A straightforward generalization of the argument in the Appendix shows that these maps are induced by the obvious inclusion maps $I^{\bar p}H_{i-r}^c(U-U\cap X^{n-k};\mc E)\to I^{\bar p}H_{i-r}^c(U;\mc E)$. Under stratum-preserving homotopy equivalence, these correspond to the inclusions $I^{\bar p}H_{i-r}^c(\mc L_x;\mc E)\to I^{\bar p}H_{i-r}^c(c\mc L_x;\mc E)$ and hence are isomorphisms for $i-r<k-1-\bar p(k)=\bar q(k)+1$ (i.e. for $i\leq \bar q(k)+r$) and $0$ otherwise. So given any echelon of constant $r+s=i$ in $E_2(U)$, either the entire echelon maps to $0$ in $E_2(U-U\cap X^{n-k})$ (when $i>\bar q(k)$), or it gets taken isomorphically to the corresponding echelon in $E_2(U-U\cap X^{n-k})$ (when $i\leq \bar q(k)$). 

Now, the trouble, of course, is that the various echelons in a spectral sequence interact as we ``turn the crank'', but, fortunately, not so badly that we can't impose some conditions that will help. We are concerned about $\H^i(U;\mc D^*(\mc I^{\bar p}\mc S^*(\mc E))[-n])\to \H^i(U-U\cap X_{n-k};\mc D^*(\mc I^{\bar p}\mc S^*(\mc E))[-n])$,  and these terms come, at the $\infty$ stage of the spectral sequence, from the echelons with $r+s=i$. We just showed that the corresponding $E_2$ terms are isomorphisms in the echeleons with $r+s=i\leq \bar q(k)$. Thus we need only impose conditions that will guarantee that these echelons continue to map isomorphically for all levels of the spectral sequence. Since each generalized boundary map $d$ at each stage of the spectral sequence never lets an echelon interact with an echelon past the one on its right, we see then that it is sufficient to force $\Ext^r(I^{\bar p}H_{\bar q(k)+1-r}^c(U-U\cap X^{n-k};\mc E),R)\cong \Ext^r(I^{\bar p}H_{\bar q(k)+1-r}^c(\mc L_x;\mc E),R)$ to be $0$ for all $r$. This corresponds to the echelon with $i=r+s=\bar q(k)+1$. Since the corresponding terms $\Ext^r(I^{\bar p}H_{\bar q(k)+1-r}^c(U;\mc E),R)\cong \Ext^r(I^{\bar p}H_{\bar q(k)+1-r}^c(c\mc L_x;\mc E),R)$ are already $0$ in $E_2(U)$, this guarantees an isomorphism at this echelon for all stages of the spectral sequence (all entries in both corresponding echelons will be $0$). This suffices to ensure then that all maps below this echelon will continue to be isomorphisms at each stage, inducing the desired isomorphism $\H^i(U;\mc D^*(\mc I^{\bar p}\mc S^*(\mc E))[-n])\to \H^i(U-U\cap X_{n-k};\mc D^*(\mc I^{\bar p}\mc S^*(\mc E))[-n])$ in the desired range.

Note that, since we consider the spectral sequence beginning at its $E_2$ stage, we in fact only need 
$\Ext^r(I^{\bar p}H_{\bar q(k)+1-r}^c(U-U\cap X^{n-k};\mc E),R)=0$ for $r\geq 2$, since no nontrivial boundary $d_j$, $j\geq 2$, will map into the $r=0$ or $r=1$ columns of the spectral sequence. This illustrates why this issue doesn't arise for principal ideal domains. Notice also that these conditions form a subset of the conditions we determined for the preceding axiom.

So, in summary, we have proven the following, which also extends the known results on pseudomanifolds:

\begin{theorem}
Let $X$ be a  MHSS with no codimension one stratum  and with sufficiently many local approximate tubular neighborhoods. Let $\mc E$ be a local coefficient system on $X-X^{n-2}$ of finitely-generated free modules over the commutative noetherian ring $R$ of finite cohomological dimension. Let $\bar p$ and $\bar q$ be dual perversities.

Suppose that for all $k$ and each $x\in X_{n-k}$, $\Ext^r(I^{\bar p}H_{i-r}^c(\mc L_x;\mc E),R)=0$  for $i>\bar q(k)$ and $r\geq 1$. Then $\mc D^*_X(\mc I^{\bar p}\mc S^*(\mc E))[-n]$ is quasi-isomorphic to $\mc I^{\bar q}\mc S^*(\mc D^*_{X-X^{n-2}}(\mc E)[-n])$.
\end{theorem}

\appendix\label{S: app}

\section{Naturality of dualization}

In this appendix, we will prove the following lemma, which is no doubt well-known but which the author has had difficulty pinpointing in the literature:

\begin{lemma}
Let $V\subset U$ be open subsets of a locally compact space $X$ of finite cohomological dimension, and let $\ms S^*\in D^b(X)$ be a sheaf complex of $R$-modules for a principal ideal domain $R$. Then there is a commutative diagram 

\begin{diagram}[LaTeXeqno]\label{naturality}
0&\rTo & \Ext(\H^{1-*}_c(U;\mc S^*),R)  &\rTo & \H^*(U;\mc D^*(\mc S^*))  &\rTo & \Hom(\H^{-*}_c(U;\mc S^*),R)   &\rTo & 0\\
 &     &\dTo & &\dTo & &\dTo &\\
0&\rTo & \Ext(\H^{1-*}_c(V;\mc S^*),R)   &\rTo & \H^*(V;\mc D^*(\mc S^*))  &\rTo & \Hom(\H^{-*}_c(V;\mc S^*),R)  &\rTo & 0,
\end{diagram}
 where the middle vertical map is induced by restriction. If $\H^*_c(W;\mc S^*)\cong H^*(\Gamma_c(W;\mc S^*))$, in particular if $X$ is an MHSS and  $\mc S^*=\mc{IS}^*$, then the side  maps are, up to isomorphism, induced by the map $\H^{-*}_c(V;\mc S^*)\to \H^{-*}_c(U;\mc S^*)$ induced by inclusion of sections. 
\end{lemma}

\begin{proof}

Recall \cite[V.7.7]{Bo} that $\mc D^*(\mc S^*)$ can be defined as the presheaf $$U\to \Hom^*(\Gamma_c(X,(\mc S^*\otimes \mc K^*)_U), I^*),$$
where $I^*$ is an injective resolution of $R$, $\mc K^*$ is an injective resolution of $R_X$ (the constant sheaf on $X$ with stalks $R$), and the subscript $U$ indicates extension by $0$ of the restriction to $U$. Note that $\mc S^*\otimes \mc K^*$ is simply a convenient $c$-soft resolution of $\mc S^*$; any $c$-soft resolution would do. So, using this definition, the restriction of sections of $\mc D^*(\mc S^*)$ from $U$ to $V$ is induced precisely by the inclusion $\Gamma_c(X,(\mc S^*\otimes \mc K^*)_V)\to \Gamma_c(X,(\mc S^*\otimes \mc K^*)_U)$. 

Now, as also noted in \cite[V.7.7]{Bo}, $\H^*(U;\mc D^*(\mc S^*))\cong \Ext^*(\Gamma_c((\mc S^*\otimes \mc K^*)_U),R)$ is the abutment of a spectral sequence with $E_2^{p,q}\cong \Ext^p(\H_c^{-q}(U;\mc S^*),R)$. This is the Cartan-Eilenberg spectral sequence for the functors $\Gamma_c(X;\cdot)$ and $\Hom^*(\cdot, I^*)$, and it is from the collapse of this spectral sequence, owing to $R$ being a principal ideal domain, that we obtain the ``universal coefficient" exact sequences that are the rows of diagram \eqref{naturality}. The restriction $(\mc S^*\otimes \mc K^*)_U\to (\mc S^*\otimes \mc K^*)_V$ induces a map of spectral sequences and hence a map of the resulting exact sequences. We just need to check that these are indeed the desired maps. But clearly the inclusion of sections $\Gamma_c(X,(\mc S^*\otimes \mc K^*)_V)\to \Gamma_c(X,(\mc S^*\otimes \mc K^*)_U)$ is 
equivalent to the inclusion $\Gamma_c(V,\mc S^*\otimes \mc K^*)\to \Gamma_c(U,\mc S^*\otimes \mc K^*)$, which gives rise to the morphism $i_*:\H^*_c(V;\mc S^*)\to \H_c^*(U;\mc S^*)$. The maps of $E_2$ terms, corresponding to the outer terms in \eqref{naturality}, are then obtained as $\Ext^*(i_*,R)$. The second claim of the lemma now follows from the natural commutative diagram
\begin{diagram}
\Gamma_c(V,\mc S^*) &\rTo & \Gamma_c(V, \mc S^*\otimes \mc K^*)\\ 
\dTo  &&\dTo\\
\Gamma_c(U,\mc S^*) &\rTo & \Gamma_c(U, \mc S^*\otimes \mc K^*).
\end{diagram}
 
For  the middle vertical map of diagram \eqref{naturality}, note that, by the naturality of the spectral sequence,  this is the morphism obtained by applying the derived functor of the composite functor $\Hom^*(\Gamma_c(X;\cdot),I^*)$ to the sheaf inclusion $(\mc S^*\otimes \mc K^*)_V\to (\mc S^*\otimes \mc K^*)_U$. But this is precisely the definition of the restriction map of the sheaf $\mc D^*(\mc S^*)$. 
\end{proof}

\bibliographystyle{amsplain}
\bibliography{bib}

\providecommand{\bysame}{\leavevmode\hbox to3em{\hrulefill}\thinspace}
\providecommand{\MR}{\relax\ifhmode\unskip\space\fi MR }
\providecommand{\MRhref}[2]{%
  \href{http://www.ams.org/mathscinet-getitem?mr=#1}{#2}
}
\providecommand{\href}[2]{#2}
\begin{thebibliography}{10}

\bibitem{Ba02}
Markus Banagl, \emph{Extending intersection homology type invariants to
  non-{W}itt spaces}, vol. 160, Memoirs of the Amer. Math. Soc., no. 760,
  American Mathematical Society, Providence, RI, 2002.

\bibitem{Ba06}
\bysame, \emph{The {L}-class of non-{W}itt spaces}, Ann. of Math. (2)
  \textbf{163} (2006), 743--766.

\bibitem{BaIH}
\bysame, \emph{Topological invariants of stratified spaces}, Springer
  Monographs in Mathematics, Springer-Verlag, New York, 2006.

\bibitem{BCS}
Markus Banagl, Sylvain Cappell, and Julius Shaneson, \emph{Computing twisted
  signatures and {L}-classes of stratified spaces}, Math. Ann. \textbf{326}
  (2003), 589--623.

\bibitem{Be97}
A.~Beshears, \emph{{G}-isovariant structure sets and stratified structure
  sets}, Ph.D. thesis, Vanderbilt University, 1997.

\bibitem{Bo}
A.~Borel~et. al., \emph{Intersection cohomology}, Progress in Mathematics,
  vol.~50, Birkhauser, Boston, 1984.

\bibitem{Br}
Glen Bredon, \emph{Sheaf theory}, Springer-Verlag, New York, 1997.

\bibitem{CS95}
Sylvain Cappell and Julius Shaneson, \emph{The mapping cone and cylinder of a
  stratified map}, Prospects in topology (Princeton, NJ), Ann. of Math. Stud.,
  vol. 138, Princeton Univ. Press, 1994, pp.~58--66.

\bibitem{CS}
Sylvain~E. Cappell and Julius~L. Shaneson, \emph{Singular spaces,
  characteristic classes, and intersection homology}, Annals of Mathematics
  \textbf{134} (1991), 325--374.

\bibitem{Cu92}
Stephen~J. Curran, \emph{Intersection homology and free group actions on {W}itt
  spaces}, Michigan Math. J. \textbf{39} (1992), 111--127.

\bibitem{GBF13}
Greg Friedman, \emph{Intersection homology of stratified fibrations and
  neighborhoods}, to appear in Advances in Mathematics; see also
  http://arxiv.org/abs/math.GT/0701112.

\bibitem{GBF3}
\bysame, \emph{Stratified fibrations and the intersection homology of the
  regular neighborhoods of bottom strata}, Topology Appl. \textbf{134} (2003),
  69--109.

\bibitem{GBF5}
\bysame, \emph{Intersection homology of regular and cylindrical neighborhoods},
  Topology Appl. \textbf{149} (2005), 97--148.

\bibitem{GBF11}
\bysame, \emph{Superperverse intersection cohomology: stratification
  (in)dependence}, Math. Z. \textbf{252} (2006), 49--70.

\bibitem{GBF10}
\bysame, \emph{Singular chain intersection homology for traditional and
  super-perversities}, Trans. Amer. Math. Soc. \textbf{359} (2007), 1977--2019.

\bibitem{GM1}
Mark Goresky and Robert MacPherson, \emph{Intersection homology theory},
  Topology \textbf{19} (1980), 135--162.

\bibitem{GM2}
\bysame, \emph{Intersection homology {II}}, Invent. Math. \textbf{72} (1983),
  77--129.

\bibitem{GS83}
Mark Goresky and Paul Siegel, \emph{Linking pairings on singular spaces},
  Comment. Math. Helvetici \textbf{58} (1983), 96--110.

\bibitem{HS91}
Nathan Habegger and Leslie Saper, \emph{Intersection cohomology of cs-spaces
  and {Z}eeman's filtration}, Invent. Math. \textbf{105} (1991), 247--272.

\bibitem{HS93}
Gilbert Hector and Martin Saralegi, \emph{Intersection cohomology of
  {$S^1$}-actions}, Trans. AMS \textbf{338} (1993), 263--288.

\bibitem{Hu99a}
Bruce Hughes, \emph{Stratifications of mapping cylinders}, Topology and Its
  Applications \textbf{94} (1999), 127--145.

\bibitem{Hug}
\bysame, \emph{Stratified path spaces and fibrations}, Proc. Roy. Soc.
  Edinburgh Sect. A \textbf{129} (1999), 351--384.

\bibitem{Hu02}
\bysame, \emph{The approximate tubular neighborhood theorem}, Ann. of Math (2)
  \textbf{156} (2002), 867--889.

\bibitem{HR}
Bruce Hughes and Andrew Ranicki, \emph{Ends of complexes}, Cambridge Tracts in
  Mathematics, vol. 123, Cambridge University Press, Cambridge, 1996.

\bibitem{HTWW}
Bruce Hughes, Laurence~R. Taylor, Shmuel Weinberger, and Bruce Williams,
  \emph{Neighborhoods in stratified spaces with two strata}, Topology
  \textbf{39} (2000), 873--919.

\bibitem{HW}
Bruce Hughes and Shmuel Weinberger, \emph{Surgery and stratified spaces},
  Surveys on Surgery Theory {V}ol. 2 (Princeton, N.J.) (Sylvain Cappell, Andrew
  Ranicki, and Jonathan Rosenberg, eds.), Annals of Mathematical Studies, vol.
  149, Princeton University Press, 2001, pp.~319--352.

\bibitem{HuW}
Witold Hurewicz and Henry Wallman, \emph{Dimension theory}, Princeton
  University Press, Princeton, 1948.

\bibitem{Ki}
Henry~C. King, \emph{Topological invariance of intersection homology without
  sheaves}, Topology Appl. \textbf{20} (1985), 149--160.

\bibitem{KIR}
Frances Kirwan, \emph{An introduction to intersection homology theory}, Pitman
  Research Notes in Mathematics Series, vol. 187, Longman Scientific and
  Technical, Harlow, 1988.

\bibitem{Bry}
Jean {L}uc Brylinski, \emph{Equivariant intersection cohmology}, Contemp. Math.
  \textbf{139} (1992), 5--32.

\bibitem{Pa05}
G.~Padilla, \emph{Intersection cohomology of stratified circle actions},
  Illinois J. Math. \textbf{49} (2005), 659--685.

\bibitem{PS04}
Gabriel Padilla, Jos\'e Ignacio~Royo Prieto, and Martintxo Saralegi-Aranguren,
  \emph{Intersection cohomology of circle actions},
  http://www.arxiv.org/abs/math.AT/0403100.

\bibitem{Q2}
Frank Quinn, \emph{Intrinsic skeleta and intersection homology of weakly
  stratified sets}, Geometry and topology (Athens, GA, 1985), Lecture Notes in
  Pure and Appl. Math., vol. 105, Dekker, New York, 1987, pp.~225--241.

\bibitem{Q1}
\bysame, \emph{Homotopically stratified sets}, J. Amer. Math. Soc. \textbf{1}
  (1988), 441--499.

\bibitem{Sa96}
Martintxo Saralegi-{A}ranguren, \emph{Cohomologie d'intersection des actions
  toriques simples}, Indag. Mathem., N.S. \textbf{7} (1996), 389--417.

\bibitem{Sa05}
Martintxo~E. Saralegi-Aranguren, \emph{De {R}ham intersection cohomology for
  general perversities}, http://www.arxiv.org/abs/math.AT/0404130.

\bibitem{Si83}
P.H. Siegel, \emph{Witt spaces: a geometric cycle theory for {KO}-homology at
  odd primes}, American J. Math. \textbf{110} (1934), 571--92.

\bibitem{Wein}
Shmuel Weinberger, \emph{The topological classification of stratified spaces},
  Chicago Lectures in Mathematics, University of Chicago Press, Chicago, IL,
  1994.

\bibitem{Ya01}
Shmuel Weinberger and Min Yan, \emph{Equivariant periodicity for abelian group
  actions}, Adv. in Geom. \textbf{1} (2001), 49--70.

\bibitem{Ya05}
\bysame, \emph{Equivariant periodicity for compact group actions}, Adv. in
  Geom. \textbf{5} (2005), 363--376.

\bibitem{Ya93}
Min Yan, \emph{The periodicity in stable equivariant surgery}, Comm. Pure Appl.
  Math \textbf{46} (1993), 1012--1040.

\end{thebibliography}

Several diagrams in this paper were typeset using the \TeX\, commutative
diagrams package by Paul Taylor. 

\end{document}